\magnification=1100 
\font\bigbf=cmbx10 scaled \magstep2
\font\medbf=cmbx10 scaled \magstep1 
\hfuzz=20 pt
\baselineskip=16pt
\input psfig.sty

\font\medbf=cmbx10 at 13pt
\def\i{\item}
\def\n{\noindent}
\def\Hat{\widehat}

\def\L{{\bf L}}
\def\meas{{\rm meas}}
\def\dint{\int\!\!\int}
\def\forall{\hbox{~~for all}~~}
\def\sqr#1#2{\vbox{\hrule height .#2pt
\hbox{\vrule width .#2pt height #1pt \kern #1pt
\vrule width .#2pt}\hrule height .#2pt }}
\def\square{\sqr74}
\def\endproof{\hphantom{MM}\hfill\llap{$\square$}\goodbreak}
\def\ve{\varepsilon}
\def\n{\noindent}
\def\c{\centerline}

\def\A{{\cal A}}

\def\vp{\varphi}

\def\O{{\cal O}}
\def\S{{\cal S}}
\def\SS{{\cal SS}}
\def\Ups{\Upsilon}
\def\BS{{\cal BS}}
\def\Rar{{\cal R}}
\def\nat{\natural}
\def\sgn{{\rm sgn}}
\def\C{{\cal C}}
\def\R{I\!\!R}

\def\wto{\rightharpoonup}
\def\osc{\hbox{Osc.}}
\def\tv{\hbox{Tot.Var.}}
\def\vs{\vskip 2em}
\def\vsk{\vskip 4em}
\def\v{\vskip 1em}
\null 
\c{\bigbf On the  Convergence Rate of}
\v
\c{\bigbf
Vanishing Viscosity Approximations} 
\vs 
\c{\it Alberto Bressan$^{(*)}$ and
Tong Yang$^{(**)}$}
\v
\c{(*) ~S.I.S.S.A., Via Beirut 4, Trieste 34014, ITALY}
\v
\c{(**)~ Department of Mathematics, City University of Hong Kong, Hong Kong}
\vsk
\n{\bf Abstract.}  Given
a strictly hyperbolic, genuinely
nonlinear system of conservation laws, we prove the a priori bound
$\big\|u(t,\cdot)-u^\ve(t,\cdot)\big\|_{\L^1}=
\O(1)(1+t)\cdot \sqrt\ve|\ln\ve|$ on the distance
between an exact BV solution $u$ and a viscous approximation
$u^\ve$, letting the viscosity coefficient $\ve\to 0$.
In the proof,
starting from $u$ we construct an approximation of the viscous 
solution $u^\ve$ by taking a mollification $u*\varphi_{\strut \sqrt\ve}$ 
and inserting
viscous shock profiles at the locations of finitely many large shocks,
for each fixed $\ve$.
Error estimates are then obtained by introducing new Lyapunov functionals 
which
control shock interactions, interactions between waves of different
families and by using
sharp decay estimates for positive nonlinear waves.

\vsk 

\n{\medbf 1 - Introduction}
\v
Consider a strictly hyperbolic system of conservation laws
$$u_t+f(u)_x=0\eqno(1.1)$$
together with the viscous approximations
$$u_t^\ve+A(u^\ve)u^\ve_x=\ve u_{xx}^\ve\,.\eqno(1.2)$$
Here $A(u)\doteq Df(u)$ is the Jacobian matrix of $f$. Given an
initial data $u(0,x)=\bar u(x)$ having small total variation, the
recent analysis in [BiB] has shown that the corresponding
solutions $u^\ve$ of (1.2) exist for all $t\geq 0$, have uniformly
small total variation and converge to a unique solution of (1.1)
as $\ve\to 0$.  The aim of the present paper is to estimate the
distance $\big\|u^\ve(t)-u(t)\big\|_{\L^1}\,$, thus
providing  a convergence rate for these vanishing viscosity
approximations. \v

We use the Landau notation $\O(1)$ to denote a quantity
whose absolute value remains uniformly bounded, while $o(1)$ indicates
a quantity that approaches zero as $\ve\to 0$.
Our main result is the following.

\v\n{\bf Theorem 1.}  {\it Let the system (1.1) be strictly
hyperbolic and assume that each characteristic field is
genuinely nonlinear. Then, given any 
initial data $u(0,\cdot)=\bar u$ with small total
variation, for every $\tau>0$
the corresponding solutions $u,u^\ve$ of (1.1) and
(1.2) satisfy the estimate}
$$\big\|u^\ve(\tau,\cdot )-u(\tau,\cdot)\big\|_{\L^1}=
\O(1)\cdot (1+ \tau)\sqrt\ve|\ln\ve|\,\tv\{\bar u\}\,.\eqno (1.3)$$ \v

\n{\bf Remark 1.} For a fixed time $\tau>0$, a similar convergence
rate was proved in [BM] for approximate solutions generated by the
Glimm scheme, namely
$$\big\|u^{Glimm}(\tau,\cdot)-u(\tau,\cdot)\big\|_{\L^1}=o(1)\cdot
\sqrt\ve\,|\ln\ve|\,.$$ Here $\ve\approx \Delta x\approx\Delta t$
measures the mesh of the grid.
\v
\n{\bf Remark 2.}
For a scalar conservation law, the method of Kuznetsov [K] shows
that the convergence rate in (1.3) is $\O(1)\cdot \ve^{1/2}$.
As shown in [TT], this rate is sharp in the general case.

In the case of hyperbolic systems, in [GX] 
Goodman and Xin have studied the viscous approximation of
piecewise smooth solutions having a finite number of
non-interacting shocks.  With these regularity assumptions,
they obtain the convergence rate $\O(1)\cdot \ve^\gamma$
for any $\gamma<1$. On the other hand, the estimate (1.3) applies
to a general BV solution, possibly with a countable everywhere dense 
set of shocks.

\vs
To appreciate the estimate in (1.3), call $S_t$ and $S_t^\ve$ the
semigroups generated by the systems (1.1) and (1.2) respectively.
As proved in [BCP], [BLY] and [BB], they are Lipschitz continuous
w.r.t.~the initial data, namely
$$\big\|S_t\bar u-S_t\bar v\big\|_{\L^1}\leq L\,
\big\|\bar u-\bar v\big\|_{\L^1}\,,\eqno(1.4)$$
$$\big\|S^\ve_t\bar u-S^\ve_t\bar v\big\|_{\L^1}\leq L\,
\big\|\bar u-\bar v\big\|_{\L^1}\,.\eqno(1.5)$$ The Lipschitz
constant $L$ here does not depend on $t,\ve$. 
By (1.4), a trivial error
estimate is
$$\big\|u^\ve(\tau)-u(\tau)\big\|_{\L^1}~=~
L\cdot\int_0^\tau\left\{\lim_{h\to 0+}{\big\|u^\ve(t+h)-S_h
u^\ve(t)\big\|_{\L^1}\over h}\right\}\,dt~=~ L\cdot\int_0^\tau
\big\|\ve u^\ve_{xx}(t)\big\|_{\L^1}\,dt\,.$$ However,
$\big\|u^\ve_{xx}(t)\big\|_{\L^1}$ grows like $\ve^{-1}$, hence
the right hand side in the above estimate does not converge to
zero as $\ve\to 0$. \v

We thus need to take a different approach, relying on (1.5).  
Let $\ve>0$ be given.
It is well known (see [B2]) that one can construct an
$\ve'$-approximate front tracking solution
$\tilde u$ of (1.1), with
$$\big\|\tilde u(0)-\bar u\big\|_{\L^1} <\ve',\qquad\qquad
\big\|\tilde u(\tau)- u(\tau)\big\|_{\L^1}
<\ve',$$ and such that the total strength of all
non-physical fronts is $<\ve'$.
Here we can take for example $\ve'=e^{-1/\ve}$.
Since the errors due to the front tracking approximation
are of order $\ve'<\!<\ve$, in the following computations we shall
neglect terms of order $\O(1)\cdot\ve'$ as   they can be 
made arbitrarily small by a suitable choice of $\ve'$.
For sake of definitiness, we shall always work with the
right-continuous version of a BV function.  Since all characteristic fields
are genuinely nonlinear, it is convenient to measure the
(signed) strength of an $i$-rarefaction or of an $i$-shock front
connecting the states $u^-,u^+$ as
$$\sigma\doteq \lambda_i(u^+)-\lambda_i(u^-)\,,$$
where $\lambda_i$ denotes the $i$-th eigenvalue of the matrix $A(u)$.
We follow here the notations in [B2], and call
$$V(u)=\sum_\alpha
|\sigma_\alpha|\,,\qquad\qquad Q(u)\doteq \sum_{(\alpha,\beta)\in\A}
|\sigma_\alpha\,\sigma_\beta|\eqno(1.6)$$
respectively the {\it total strength of waves} and the 
{\it interaction potential} in a front tracking
solution $u$.   The second summation
here ranges over the set $\A$ of all couples of approaching wave fronts.

For notational convenience,
we shall simply call $u$ the $\ve'$-approximate
front tracking approximation, also
assume that $\bar u=u(0)$ is piecewise constant. Since
$\ve'<\!<\ve$, this will not have any consequence for our estimates. 
In the sequel, we shall construct a further
approximation $v=v(t,x)$ having the following properties. 
\v
Let $0=t_0<t_1<\cdots<t_N=\tau$ be the interaction times in the
front tracking solution $u$.    Then $v$ is smooth
on each strip
$[t_{i-1}\,,~t_i[\,\times \R$.  Moreover, calling $\delta_0\doteq
\tv\{\bar u\}$, one has
$$\big\|v(0)-\bar u\big\|_{\L^1}=\O(1)\cdot\delta_0 \sqrt\ve\,,\qquad\qquad
\big\|v(\tau)-u(\tau)\big\|_{\L^1}=\O(1)\cdot\delta_0
\sqrt\ve\,,
\eqno(1.7)
$$
$$\int_0^\tau \int \big| v_t+A(v)v_x-\ve v_{xx}\big|\,dxdt
=\O(1)\cdot \delta_0
(1+\tau)\sqrt\ve\,|\ln\ve|\,,
\eqno(1.8)
$$
$$\sum_{1\leq i\leq N} \int \big|v(t_i,x)-
v(t_i-,x)\big|\,dx=\O(1)\cdot \delta_0\sqrt\ve\,|\ln\ve|\,.
\eqno(1.9)
$$
Having achieved this step, by the Lipschitz continuity of the
semigroup $S_t^\ve$ in (1.5) we can then conclude
$$\eqalign{\big\|u^\ve(\tau)-u(\tau)\big\|_{\L^1}&\leq
\big\|S^\ve_\tau \bar
u-v(\tau)\big\|_{\L^1}+\big\|v(\tau)-u(\tau)\big\|_{\L^1}\cr &\leq
L\,\big\|\bar u-v(0)\big\|_{\L^1}+L\, \int_0^\tau\!\int \big|
v_t+A(v)v_x-\ve v_{xx}\big|\,dxdt\cr &\qquad +L\,\sum_{1\leq i\leq
N}
\big\|v(t_i,x)-v(t_i-,x)\big\|_{\L^1}+\big\|v(\tau)-u(\tau)\big\|_{\L^1}
\cr &=\O(1)\cdot \delta_0(1+\tau)\sqrt\ve|\ln\ve|\,.\cr}\eqno(1.10)$$

To construct the approximate solution $v$, we first consider a
mollification of $u$ w.r.t.~the space variable $x$. Let
$\vp:\R\mapsto [0,1]$ be a smooth function such that
$$\vp(s)=0\qquad\hbox{if}~~|s|>{2\over 3}\,\qquad \qquad s\,\vp'(s)\leq
0,\qquad \vp(s)=\vp(-s),\quad\qquad\int\vp(s)\,ds=1\,.$$ For
$\delta>0$ small, define the rescalings $\vp_\delta(s)\doteq
\delta^{-1}\vp(x/\delta)$ and the mollified solutions
$v^\delta(t)\doteq u(t)*\vp_\delta$, so that
$$v^\delta(t,x)=\int u(t,y)\vp_\delta(x-y)\,dy\,.$$
Recalling that $\delta_0\doteq\tv\{\bar u\}$, one has
$$\tv\big\{ u(t)\big\}\,,~~\big\|u^\ve_x(t)\big\|_{\L^1}\,=
\,\O(1)\cdot\delta_0,\qquad\quad\forall t\geq 0\,. \eqno(1.11)$$ We
now observe that
$$\eqalign{\big\|v^\delta-u\big\|_{\L^1}&=\int\left|
\int\big(u(x)-u(y)\big)\vp_\delta(x-y)\,dy\right|\,dx \cr
&\leq\int \tv\big\{ u\,;~[x-\delta,\,x+\delta]\big\}\,dx~=~\O(1)\cdot
\delta_0\, \delta\,.\cr}\eqno(1.12)$$ To estimate the distance
between $v^\delta$ and $u^\ve$, we first compute
$$\int\big|\ve
\,v^\delta_{xx}(x)\big|\,dx=\ve\,\int\big|(u_x*\vp_{\delta,x})(
x)\big|\,dx\leq\ve\|u_x\|_{\L^1}
\cdot\big\|\vp_{\delta,x}\big\|_{\L^1}=\O(1)\cdot \delta_0\,
{\ve\over\delta}\,.\eqno(1.13)
$$
$$
\eqalign{\int \big|v^\delta_t+A(v^\delta)v^\delta_x\big|\,dx
&=\int\left|\int
\Big(A\big(v^\delta(x)\big)u_x(y)-A\big(u(y)\big)u_x(y)\Big)
\vp(x-y)\,dy\right|\,dx
\cr
&\leq\int\left(\int\Big|A\big(v^\delta(x)\big)-A\big(u(y)\big)
\Big|\vp(x-y)\,dx\right)
\big|u_x(y)\big|\,dy\cr &=\O(1)\cdot\|DA\|_{\C^0}\int\osc\big\{
u\,;~~[y-\delta,~y+\delta]\big\}\,\big|u_x(y)\big|\,dy\,.
\cr}
\eqno(1.14)
$$
For simplicity, the formulas (1.13)-(1.14)
are here written in the case where the function $u$ is absolutely
continuous.  In the general case, the same 
estimates hold, by replacing $|u_x|dx$ with the measure
$|D_x u|$ of total variation of $u\in BV$.
\v
If $u$ is a Lipschitz continuous solution of (1.1), the
oscillation of $u$ on any interval of length $2\delta$ is
$\O(1)\cdot\delta$. Hence, performing the above mollifications,
we would obtain
$$\int\big|v^\delta_t+A(v^\delta)v^\delta_x\big|\,dx=\O(1)\cdot
\delta\,\delta_0\,.\eqno(1.15)$$ Choosing $\delta\doteq \sqrt\ve$,
by (1.12)--(1.15) we thus conclude
$$\eqalign{\big\|u^\ve(\tau)-u(\tau)\big\|_{\L^1}&\leq
\big\|S^\ve_\tau \bar
u-v^\delta(\tau)\big\|_{\L^1}+\big\|v^\delta(\tau)-u(\tau)\big\|_{\L^1}\cr
&\leq L\,\big\|\bar u-v^\delta(0)\big\|_{\L^1}+L\,
\int_0^\tau\!\int \big| v^\delta_t+A(v^\delta)v^\delta_x-\ve
v^\delta_{xx}\big|\,dxdt+\big\|v^\delta(\tau)-u(\tau)\big\|_{\L^1} \cr
&=\O(1)\cdot \delta_0(1+\tau)\sqrt\ve\,.\cr}
\eqno(1.16)$$

In general, however, the solution $u$ is not Lipschitz continuous.
The best one can say is that $u$ is 
a function with bounded variation, possibly with countably
many shocks. Hence the easy estimate (1.16) does not hold. For
genuinely nonlinear systems, the additional error terms due to
centered rarefaction waves can be controlled by carefully
estimating the decay rate of these waves. Error terms due to small
shocks will be estimated by suitable Lyapunov functionals.
However, there is one type of wave-fronts which is responsible for
large errors in (1.14), namely the large shocks of strength
$>\!>\sqrt\ve$. In a neighborhood of each one of these shocks, a
more careful approximation is needed. Instead of a mollification,
we shall insert an approximate viscous shock profile.

Our construction goes as follows. By the same argument as in [BC1]
(see Proposition 2 on p.17), given
$\rho>0$ one can select a finitely many shock fronts
$$t\mapsto x_\alpha(t)\qquad\qquad t\in
T_\alpha\doteq \,]t_\alpha^-,~t_\alpha^+[\,,\qquad\qquad
\alpha=1,\ldots,\nu,$$ with $\nu=\O(1)\cdot \delta_0/\rho$, having
the following properties.
\v
\i{$\bullet$} For every $t\in T_\alpha$ (apart from
finitely many interaction points) the left and right states
$u^-_\alpha, u^+_\alpha$ are connected by a shock, say of the
family $k_\alpha$, with strength $\big|\sigma_\alpha(t)\big|\geq \rho/2$,
while $\big|\sigma_\alpha(t^*)\big|\geq \rho$ for
some $t^*\in T_\alpha$.  Moreover,
every shock in the front tracking solution $u$ with
strength $\geq \rho$ is included in one of the above fronts. 

\v
For each $\alpha$ and each $t\in T_\alpha$ (apart from
finitely many interaction points), let $\omega_\alpha$ be
the viscous shock profile connecting the states
$u_\alpha^-,u_\alpha^+$.   Calling $\lambda_\alpha$ the shock
speed, we thus have
$$
\omega''_\alpha=\big(A(\omega_\alpha)-\lambda_\alpha\big)
\omega'_\alpha\,,\qquad\qquad
\lim_{s\to\pm\infty}\omega_\alpha(s)=u_\alpha^\pm\,.
$$
We choose the parameter $s$ so that the value $s=0$
corresponds roughly to the center of the travelling profile.  This
can be achieved by requiring
$$\int_{-\infty}^0\big|\omega_\alpha(s)-u_\alpha^-\big|\,ds=
\int_0^\infty\big|\omega_\alpha(s)-u_\alpha^+\big|\,ds\,.
\eqno(1.17)
$$
For the system (1.2) with $\ve$-viscosity, the
corresponding rescaled shock profile is
$s\mapsto\omega_\alpha^\ve(s)\doteq \omega_\alpha(s/\ve)$. On the
open interval $$J_\alpha(t)\doteq
~\big]x_\alpha(t)-\delta\,~,~~x_\alpha(t)+\delta\big[$$ we now
replace the mollified solution by a shock profile. Define the
functions $\varrho_\alpha$, $\tilde\omega_\alpha$, by setting
$$\varrho_\alpha(x_\alpha+\xi)\doteq
u_\alpha^+\int_{-\infty}^\xi\vp_\delta(y)\,dy+
u_\alpha^-\int_\xi^\infty\vp_\delta(y)\,dy\, ,
\eqno(1.18)$$
$$
\tilde\omega_\alpha(x_\alpha+\xi)\doteq\cases{\omega^\ve_\alpha
\big( \phi(\xi))\qquad &if\quad
$\xi\in\,]-\delta\,,~\delta[\,$,\cr u_\alpha^+\qquad &if\quad
$\xi\geq \delta$,\cr u_\alpha^-\qquad &if\quad $\xi\leq
-\delta$,\cr}
\eqno(1.19)$$
where
$$
\phi(\xi)=\cases{\xi\qquad &if\quad $|\xi|\leq 
{{\sqrt{\ve}}\over{2}}\, $,\cr
&\cr
{\ve\over 4(\sqrt \ve -\xi)} \qquad &if\quad
${\sqrt\ve\over 2}\le \xi<\sqrt\ve\,$,\cr &\cr
-{\ve\over 4(\sqrt \ve +\xi)} \qquad &if\quad
$-\sqrt\ve<\xi<-{\sqrt\ve\over 2}\,$.
\cr}
\eqno(1.20)$$
Notice that
$\tilde\omega_\alpha$ is essentially an $\ve$-viscous shock
profile, up to a $\C^1$ tranformation that squeezes the whole real line
onto the interval $J_\alpha(t)$.  Moreover, $\varrho_\alpha$ is the
mollification of the piecewise constant function taking values
$u_\alpha^-,u_\alpha^+$ with a single jump at $x_\alpha$. The
above definitions imply that $\tilde\omega_\alpha=\varrho_\alpha$
outside the interval $J_\alpha(t)$. Finally, for every $t\geq 0$ we
define
$$v=u*\vp_\delta~+\,\sum_{\alpha\in\BS}
\big(\tilde\omega_\alpha-\varrho_\alpha)\,,\eqno(1.21)$$ 
where the summation ranges over all big shock fronts.
In the
remainder of the paper we will show that, by choosing
$$\delta\doteq\sqrt \ve\,,\qquad\qquad \rho\doteq
4\sqrt\ve\,|\ln\ve|\,,\eqno(1.22)$$ all the estimates in
(1.7)--(1.9) hold. By (1.10), this will achieve a proof of Theorem
1. 
\vsk

\n{\medbf 2 - Estimates on rarefaction waves} 
\v
Throughout the following we denote by 
$\lambda_1(u)<\cdots<\lambda_n(u)$ the eigenvalues
of the $A(u)\doteq Df(u)$. 
Moreover, we shall use bases of left and right eigenvectors $l_i(u)$, $r_i(u)$
normalized so that
$$\nabla\lambda_i(u)\cdot r_i(u)\equiv 1\,,\qquad\qquad 
l_i(u)\cdot r_j(u)=\cases{1\quad &if\quad $i=j$,\cr
0\quad &if\quad $i\not=j$.\cr}\eqno(2.1)$$
According to (1.14), outside the large shocks we have to estimate
the quantity
$$E(\tau)\doteq \int_0^\tau\!\int\osc\big\{
u\,;~~[y-\delta,~y+\delta]\big\}\,\big|u_x(y)\big|\,dydt\,.\eqno(2.2)$$
Centered rarefaction waves can have large gradients, and hence
give a large contribution to the above integral. However, for
genuinely nonlinear families, the density of these waves decays
rapidly, as $t^{-1}$. We now give an example where 
the integral (2.2) can be easily estimated.
\v
\n{\bf Example 1.}  Assume that the solution $u$ consists of a single
centered rarefaction wave of the $i$-th family (fig.~1), connecting the
states $u^-,u^+$.
Call $s\mapsto \omega(s)$ the parametrized $i$-rarefaction curve,
so that
$$\dot\omega=r_i(\omega),\qquad\qquad\omega(0)=u^-\,,\quad
\omega(\sigma)=u^+$$
for some wave strength $\sigma>0$.
We then have
$$u(t,x)=\cases{ u^-\qquad &if\quad $x/t<\lambda_i(u^-)\,,$ \cr
\omega(s)\qquad &if\quad $x/t=\lambda_i\big(\omega(s)\big)\,\qquad
s\in [0,\sigma]$,\cr
u^+\qquad &if\quad $x/t>\lambda_i(u^+)$.\cr}$$
If $K$ is an upper bound for the length of all eigenvectors $r_i(u)$,
we have
$$\osc\big\{u(t) \,;~~[y-\delta,~y+\delta]\big\}\leq K\cdot\min\big\{
\sigma\,,~2\delta/t\big\}\,,\qquad\quad \int \big|u_x(x)\big|\,dx\leq K\sigma
\,.$$ 
Hence the quantity in (2.2) satisfies
$$E(\tau)~\leq~\int_0^{2\delta/\sigma}K^2
\sigma^2\,dt+\int_{2\delta/\sigma}^\tau
K^2\sigma{2\delta\over t}\,dt~=~2K^2\,\delta\sigma\left(1+
\ln{\sigma\tau\over 2\delta}\right)\,.\eqno(2.3)$$ 
The
choice $\delta\doteq \sqrt\ve$ would thus give the correct
order of magnitude $\O(1)\cdot \sqrt\ve\,|\ln\ve|\cdot
\tv\{u\}$.

\midinsert
\vskip 10pt
\centerline{\hbox{\psfig{figure=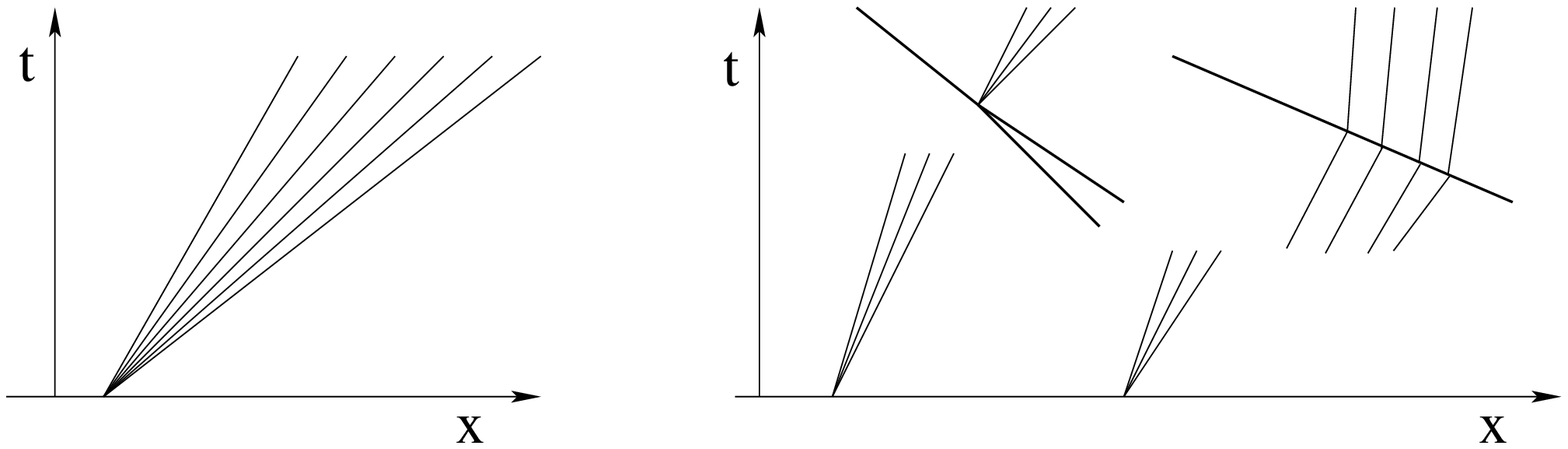,width=14cm}}}
\centerline{\hbox{figure 1~~~~~~~~~~~~~~~~~~~~~
~~~~~~~~~~~~~~~~~~~~~~~~~~~~~figure 2~~~~~~~~~~~~}}
\vskip 10pt
\endinsert

Of course, a general BV solution of the system of conservation laws (1.1)
is far more complex than a single rarefaction.  It can contain
several centered rarefactions originating at $t=0$ and also
at later times, as a result of shock interactions (fig.~2).
Moreover, the crossing of wave fronts of other families may
slow down the decay of positive waves.   Nevertheless,
the forthcoming analysis will show that, in some sense, 
Example 1 represents the worst possible case.   
Using the sharp decay
estimate for positive waves in [BY] and a comparison argument, we
shall prove that the total error due to steep rarefaction waves
for an arbitrary weak solution is no greater than the error
computed at (2.3) for a solution containing only one centered
rarefaction. In the present section, all the analysis refers to
an exact solution.  A similar result can then be easily derived
for a sufficiently accurate front tracking approximation.

We begin by recalling the main results in [BY].
Given a function $u:\R\mapsto\R^n$ with small total variation, 
following [BC] and [B2], one can define the measures
$\mu^i$ of $i$-waves in $u$ as follows.
Since $u\in BV$, its distributional derivative $D_xu$ is a Radon measure.
We define $\mu^i$ as the measure such that
$$\mu^i\doteq l_i(u)\cdot D_xu\eqno(2.4)$$
restricted to the set where $u$ is continuous, while, at each point $x$ 
where $u$ has a jump, we define
$$\mu^i\big(\{x\}\big)\doteq\sigma_i\,,\eqno(2.5)$$
where $\sigma_i$ is the strength of the $i$-wave in the solution of
the Riemann problem with data 
$u^-=u(x-)$, $u^+=u(x+)$. 
In accordance with (2.1), if the solution of the Riemann problem contains 
the intermediate states $u^-=\omega_0, \omega_1,\ldots,\omega_n=u^+$, 
the strength of the $i$-wave is defined as 
$$\sigma_i\doteq \lambda_i(\omega_i)-\lambda_i(\omega_{i-1}).\eqno(2.6)$$ 
Together with the measures $\mu^i$ we also define the 
Glimm functionals 
$$V(u)\doteq \sum_i |\mu^i|(\R)\,,$$ 
$$Q(u)\doteq \sum_{i<j} \big(|\mu^j|\otimes|\mu^i|\big) \big\{ 
(x,y)\,;~x<y\big\} 
+\sum_i\big(\mu^{i-}\otimes|\mu^i|\big) \big\{ 
(x,y)\,;~x\not=y\big\}\,,
$$ 
measuring respectively
the total strength of waves and the interaction potential.

We call $\mu^{i+}$, $\mu^{i-}$ respectively  the positive  
and negative parts of $\mu^i$, so that 
$$\mu^i=\mu^{i+}-\mu^{i-},\qquad\qquad |\mu^i|=\mu^{i+}+\mu^{i-}.
\eqno(2.7)$$ 
In [BY], the authors
introduced a partial ordering within the  
family of positive Radon measures:
\v 
\n{\bf Definition 1.} {\it  
Let $\mu,\mu'$ be two positive Radon measures. 
We say that 
$\mu\preceq \mu'$ if and only if  
$$\sup_{meas(A)\leq s}\mu(A)~\leq~\sup_{meas(B)\leq s}\mu'(B) 
\qquad\qquad\hbox{for every}~s>0\,. 
\eqno(2.8)$$ 
}
\v 
Here meas$(A)$ denotes the 
Lebesgue measure of a set $A$. 
In some  sense, the 
above relation means that $\mu'$ is more singular than $\mu$.  
Namely, it has a greater total mass, concentrated   
on regions with higher density. 
Notice that the usual order relation 
$$\mu\leq\mu'\qquad\hbox{if and only if}\qquad \mu(A)\leq \mu'(A) \quad 
\hbox{for every}\quad A\subset\R$$ 
is much stronger. Of course 
$\mu\leq\mu'$ implies $\mu\preceq\mu'$, 
but the converse does not hold. 
 
Given a solution $u$ of (1.1), we denote by  $\mu^{i+}_t$ the
measure of positive $i$-waves in  
$u(t,\cdot)$.
In particular, $ \mu^{i+}_0$ refers to the positive $i$-waves in  
$u$ at the initial time $t=0$.  An accurate estimate of these measures
is obtained by a comparison with a solution 
of Burgers' equation with source terms. 
\v 
\n{\bf Proposition 1.} {\it  For some constant 
$\kappa>0$ and for every small BV solution $u=u(t,x)$ 
of the system (1.1) the following holds. 
Let $w=w(t,x)$ be the solution of the  Cauchy problem for
Burgers' equation 
with impulsive source term 
$$w_t+(w^2/2)_x=-\kappa\,\sgn(x)\cdot {d\over dt}Q\big(u(t)\big)\,,
\eqno(2.9)$$ 
$$w(0,x)=\sgn(x)\cdot \sup_{meas(A)<2|x|} 
{\mu^{i+}_0(A)\over 2}\,.\eqno(2.10)$$ 
Then, for every $t\geq 0$,  
$$\mu^{i+}_t\preceq D_x w(t)\,.\eqno(2.11)$$ 
} 
\v 
For a proof, see [BY].

The ordering relation (2.8) 
can be better appreciated in terms of 
rearrangements. 
More precisely,
let $\mu$ be a positive Radon measure on $\R$, so that 
$\mu \doteq D_x v$ is the distributional derivative of some 
bounded, non-decreasing function 
$v:\R\mapsto\R$. 
We can decompose 
$$\mu=\mu^{\rm sing}+\mu^{ac}$$ 
as the sum of a singular and an absolutely continuous part, w.r.t.~Lebesgue 
measure.  The absolutely continuous part corresponds to 
the usual derivative $z\doteq v_x$, which is a non-negative $\L^1$ function 
defined at a.e.~point. We shall denote by $\hat z$ the {\it 
symmetric rearrangement} of $z$, i.e.~the unique even function 
such that 
$$\hat z(x)=\hat z(-x)\,,\qquad\qquad \hat z(x)\geq \hat z(x')~~~ 
\hbox{if}~~~0<x<x'\,, 
$$ 
$$ 
\meas\Big(\big\{x\,;~\hat z(x)> c\big\}\Big)= 
\meas\Big(\big\{x\,;~z(x)> c\big\}\Big),\qquad\qquad \hbox{for every}~ c>0 
\,. 
$$ 
Moreover, we define the {\it odd rearrangement} of $v$ as the  
unique function 
$\hat v$ such that 
$$\hat v(-x)=-\hat v(x)\,,\qquad\qquad 
\hat v(0+) = {1\over 2}\mu^{\rm sing}(\R)\,,$$ 
$$\hat v(x)=\hat v(0+)+\int_0^x z(y)\,dy \qquad 
\hbox{for}~x>0\,.$$ 
By construction, the function $\hat v$ is convex for $x<0$ and  
concave for $x>0$.  We now have 
 \v 
\n{\bf Proposition 2.} {\it Let $\mu=D_x v$ and 
$\mu'=D_x v'$ be  positive Radon measures. 
Call $\hat v,\hat v'$ the odd rearrangements
of $v,v'$, respectively.   Then $\mu\preceq D_x\hat v\preceq\mu$.
Moreover
$$\hat v(x)=\sgn(x)\cdot\sup_{meas (A)\leq 2|x|} 
{\mu(A)\over 2}\,.
\eqno(2.12)$$ 
Moreover,
$$\mu\preceq\mu'\qquad\hbox{if and only if}\qquad
\hat v(x)\leq \hat v'(x)\quad\hbox{for all}~x>0\,.\eqno(2.13)$$
} 

The relevance of the above concepts toward an estimate of the quantity in
(2.2) is due to the next three comparison lemmas.
\v
\n{\bf Lemma 1.} {\it Let $u:\R\mapsto\R$ be a non-decreasing BV
function and let $\hat u$ be its odd rearrangement. Then}
$$\int_{-\infty}^\infty
\tv\big\{u\,;~[x-\rho,\,x+\rho]\big\}\,du(x)\leq
3\int_{-\infty}^\infty \big[\hat u(x+\rho)-\hat u(x-\rho)\big]
\,d\hat u(x)\,.\eqno(2.14)$$ 
\v
\n{\bf Proof.} 
We begin by defining a measurable
map $x \mapsto \vp(x)$ from $\R$ onto $\R_+$
with the following properties.
\v
\i{(i)} $\vp(x)=0$ for all points $x$ in the support of
singular part of the measure $u_x$.
\v
\i{(ii)} $u_x(x)=\hat u_x\big(\vp(x)\big)$ for every $x$ where $u$ is 
differentiable.
\v
\i{(iii)} $\meas\big(\vp^{-1}(A)\big)=2\,\meas(A)$ for every $A\subset\R_+$.
\v
We now have
$$\eqalign{ &\int_{-\infty}^\infty
\tv\big\{u\,;~[x-\rho,\,x+\rho]\big\}\,du(x)\cr
&=\left(\int_{\vp(x)\le \rho}
+\int_{\vp(x)>\rho}\right)
\big[u(x+\rho)-u(x+\rho)\big]\,du(x)\cr
&\doteq I_1+I_2.}
$$
We now estimate $I_1$ and $I_2$ separately as follows.
$$\eqalign{I_1&=\int_{\vp(x)\le\rho}
\big[u(x+\rho)-u(x-\rho)\big]\,du(x)\cr
&\leq \int_{\vp(x)\le\rho}
2\hat u(\rho)\,du(x)\cr
&\leq 4\big(\hat u(\rho)\big)^2\cr
&\leq 2\int_{-\rho}^\rho
\big[\hat u(x+\rho)-\hat u(x-\rho)\big]\,d\hat u(x).\cr}
\eqno(2.15)$$
$$\eqalign{I_2&\le\int_{\vp(x)>\rho}\int_{-\rho}^\rho
\big[u_x(x)\,D_x u(x+s)\big]\,dsdx\cr
&\leq  4\rho\int_{\rho}^\infty \big[\hat u_x(x) \,D_x\hat u(x-\rho)\big]
\,dx\cr
&= 4\rho\int_0^\infty \hat u_x(x+\rho)\,d\hat u(x)\cr
&\leq 2\int_0^\infty \big[\hat u(x+\rho)-\hat u(x-\rho)\big]
\,d\hat u(x)\cr& =
\int_{-\infty}^\infty
\big[\hat u(x+\rho)-\hat u(x-\rho)\big]\,d\hat u(x).\cr}
\eqno(2.16)
$$
For $x>\rho$, we are here using the inequality
$$2\rho \hat u_x(x)\leq \hat u(x)-\hat u(x-2\rho).$$
Moreover, calling $\tilde f,\tilde g$ the
non-increasing even rearrangements of two positive, integrable 
functions $f,g$,
one always has
$$\int_{-\infty}^\infty 
f(x)\,g(x)\,dx\leq\int_{-\infty}^\infty \tilde f(x)\,\tilde g(x)\,dx\,.
\eqno(2.17)$$
Together, (2.15) and (2.16) yield (2.14).
\endproof
\v
\n{\bf Lemma 2.} {\it 
Let $v,w$ be two non-decreasing BV functions. If $D_xv\preceq D_x w$
then the odd rearrangements $\hat v,\hat w$ satisfy}
$$
\int_{-\infty}^\infty \big[\hat v(x+\rho)-\hat v(x-\rho)\big]
\,d\hat v(x) \leq 
\int_{-\infty}^\infty\big[\hat w(x+\rho)-\hat w(x-\rho)\big]
\,d\hat w(x)\,.\eqno(2.18)$$
\v
\n{\bf Proof.} By an approximation argument, we can assume that
$\hat v$ and $ \hat w$ are smooth. Without loss of generality,
we can assume $\hat{v}(\pm\infty)=\hat{w}(\pm\infty)$.
By assumptions, 
$\hat v(x) \leq \hat w(x)$  for all $x>0$.  We consider a 
parabolic equation with smooth
coefficients
$$z_t=a(t,x)z_{xx}\,,\eqno(2.19)$$
with 
$a(t,x)=a(t,-x)\geq  0$,
having a solution such that
$$z(0,x)=\hat w(x),\qquad 
\lim_{t\to\infty}z(t,x)=\hat v(x)\,,
$$
where the limit holds uniformly for $x$ in bounded sets. 
To construct $a(t,x)$, one can first define a smooth function
$\tilde a=\tilde a(t,x,z)$ such that
$$\tilde a(t,\,-x,z)=\tilde a(t,x,z)
=\cases{1\qquad &if\quad $\big|z-\hat v(x)\big|\geq 2/t\,$,\cr
0\qquad &if\quad  $\big|z-\hat v(x)\big|\leq 1/t\,$.\cr}$$
Then we solve the quasilinear Cauchy problem
$$z_t=\tilde a(t,x,z)z_{xx}\,,\qquad\qquad z(0,x)=\hat w(x)$$
and set $a(t,x)\doteq \tilde a\big(t,x,z(t,x)\big)$.
We now claim that
$${d\over dt} \left(\int_{-\infty}^\infty\int_{x-\rho}^{x+\rho}z_x(x)\,z_x(y)
\,dydx\right)\leq 0\,.\eqno(2.20)$$
Indeed, calling $\phi\doteq z_x\geq 0$ and using (2.19) we compute
$$\phi_t=\big(a(t,x)\phi_x\big)_x\,,$$
$$\eqalign{{d\over dt} &\left(\int_{-\infty}^\infty  
\int_{x-\rho}^{x+\rho}\phi(x)\,\phi(y)
\,dydx\right)= \int_{-\infty}^\infty \int_{x-\rho}^{x+\rho}\Big[
\big(a\, \phi_x(x)\big)_x \phi(y)+\phi(x) (a\,\phi_x(y)\big)_x\Big]\,dydx
\cr
&=\int_{-\infty}^\infty \big[ a\,\phi_x(y+\rho)
-a\phi_x(y-\rho)\big] \phi(y)\,dy
+\int_{-\infty}^\infty 
\big[ a\,\phi_x(x+\rho)-a\phi_x(x-\rho)\big] \phi(x)\,dx
\cr
&=2\int_{-\infty}^\infty \phi(x)\big[a \phi_x (x+\rho)-a\phi_x(x-\rho)\big]
\,dx\cr
&=\int_{-\infty}^\infty a \phi_x(x)\big[\phi(x-\rho)-\phi(x+\rho)\big]\,dx
\cr
&\leq 0\,,\cr}$$
because $\phi(t,\cdot)$ is an even function, non-increasing for $x\geq 0$. 
{}From (2.20) it follows
$$
\int_{-\infty}^\infty\int_{x-\rho}^{x+\rho}\hat v_x(x)\,\hat v_x(y)
\,dydx\leq 
\int_{-\infty}^\infty\int_{x-\rho}^{x+\rho}\hat w_x(x)\,\hat w_x(y)
\,dydx\,.$$
\endproof
\v
\n{\bf Lemma 3.} {\it Let $u$ be a solution of (1.1) defined 
for $t\in [0,\tau]$ and let
$w=w(t,x)$ as in (2.9)-(2.10).
Set
$$\bar\sigma\doteq {1\over 2}\mu^{i+}_0(\R)+
\kappa \big[Q(u(0))-Q(u(\tau))\big]\eqno(2.21)$$
and let
$$v(t,x)= \cases{ x/t\qquad &if\qquad $|x|/t\leq \bar\sigma\,$,\cr
\sgn(x)\cdot \bar\sigma &if\qquad $|x|/t> \bar\sigma\,$,\cr}
\eqno(2.22)$$ 
be a solution of Burgers' equation consisting of one single
centered rarefaction wave of strength $2\bar\sigma$.
Then
$$
\int_0^\tau\int_{-\infty}^\infty \big[w(t\,,~x+\rho)-w(t\,,~x-\rho)\big]
\,w_x(t,x)\,dxdt \leq 2\int_0^\tau
\int_{-\infty}^\infty\big[v(t\,,~x+\rho)-v(t\,,~x-\rho)\big]v_x(t,x)
\,dx\,dt\,.\eqno(2.23)$$
}
\v
\n{\bf Proof.}
To compare the integrals in (2.23)
a change of variables will be useful.
We define (fig.3)
$$x(t,\xi)\doteq t\xi\qquad\qquad \xi\in [0,\bar\sigma]\,,\quad
t\in [0,\tau]\,.$$
For $t\in [0,\tau]$ and $ Q(t)-Q(\tau)<\xi\leq\bar\sigma$,
we also consider the point $y(t,\xi)>0$ implicitly defined by
$$w(t,\infty)-w\big(t,\,y(t,\xi)\big)=\bar\sigma-\xi\,.$$
Notice that $y(t,\xi)$ is defined only for $t\in\big[ t(\xi)\,,~\tau
\big]$,
or equivalently $\xi\in \big[\xi(t)\,,~\bar\sigma\big]$, where
$$\xi(t)\doteq \kappa\big[ Q(t)-Q(\tau)\big]\,,\qquad\qquad
t(\xi)\doteq \inf\Big\{ t\geq 0\,;~~\big[Q(t)-Q(\tau)\big]\leq\xi\Big\}\,.$$
For $0<\xi_1<\xi_2<\bar\sigma$ and $s>0$ 
we have
$$\eqalign{y\big(t(\xi_1)+s\,,~\xi_2\big)-y\big(t(\xi_1)+s\,,~\xi_1\big)
&=y\big(t(\xi_1)\,,~\xi_2\big)-y\big(t(\xi_1)\,,~\xi_1\big)+(\xi_2-\xi_1)s\cr
&\geq (\xi_2-\xi_1)s\cr
&=x(s,\xi_2)-x(s,\xi_1)\,.\cr}\eqno(2.24)$$
Observe that, since $w$ is odd and non-decreasing,
$$w^+(t,\,y-\rho)\doteq\max\big\{ w(t,\,y-\rho)\,,~0\big\}
=w\big( t\,,~\max\{ y-\rho\,,~0\}\big)\,.$$
Of course, the same is  also true for $v$.
Calling $I_w$, $I_v$ the two integrals in (2.23) and using
(2.24) at the key step,
we obtain
$$\eqalign{I_w&=2\int_0^\tau \int_{\xi(t)}^{\bar\sigma}
\Big[w\big(t\,,~y(t,\xi)+\rho\big)-w\big(t\,,~y(t,\xi)-\rho\big)
\Big]\,d\xi\,dt \cr
&\leq 4\int_0^\tau \int_{\xi(t)}^{\bar\sigma}
\Big[w\big(t\,,~y(t,\xi)+\rho\big)-w^+\big(t\,,~y(t,\xi)-\rho\big)
\Big]\,d\xi\,dt 
\cr
&=4\int_0^\tau \dint_{\big|y(t,\xi_1)-y(t,\xi_2)\big|<\rho}
d\xi_1 d\xi_2\,dt\cr
&=4\dint \hbox{meas}\,\Big\{t\in [0,\tau]\,;~~
\big|y(t,\xi_1)-y(t,\xi_2)\big|<\rho\Big\}\,d\xi_1 d\xi_2\cr
&\leq 4\dint \hbox{meas}\,\Big\{t\in [0,\tau]\,;~~
\big|x(t,\xi_1)-x(t,\xi_2)\big|<\rho\Big\}\,d\xi_1 d\xi_2\cr
&=
4\int_0^\tau \int_0^{\bar\sigma}
\Big[v\big(t\,,~x(t,\xi)+\rho\big)-v^+\big(t\,,~x(t,\xi)-\rho\big)
\Big]\,d\xi\,dt \cr
&\leq 2I_v\,.\cr}
$$
\endproof

\midinsert
\vskip 10pt
\centerline{\hbox{\psfig{figure=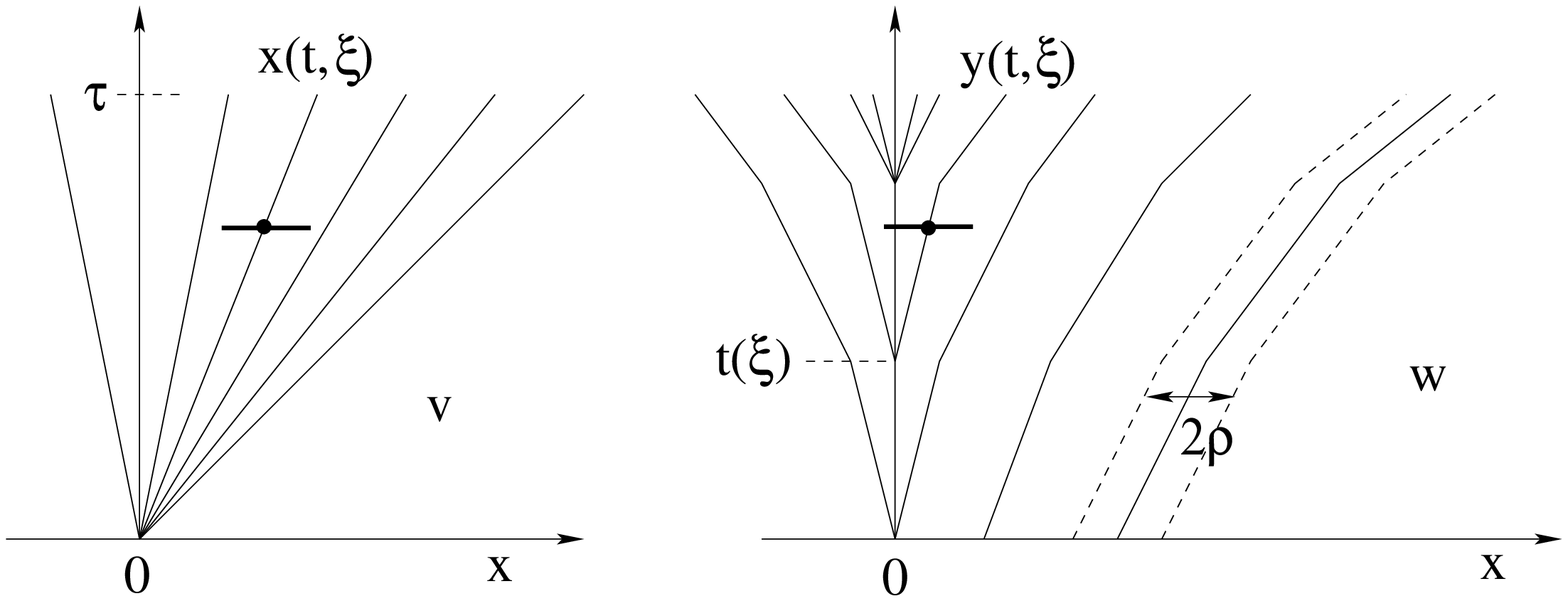,width=14cm}}}
\centerline{\hbox{figure 3}}
\vskip 10pt
\endinsert

\n{\bf Corollary 1.} {\it Assume that all 
characteristic fields for the system (1.1) are genuinely nonlinear.
Let $u$ be a solution with initial data $u(0,x)=\bar u(x)$
having small total variation.  Then, for every $\tau,\delta>0$, 
the measures $\mu^{i+}_t$ of positive
waves in $u(t,\cdot)$ satisfy the estimate}
$$\sum_{i=1}^n\int_0^\tau (\mu_t^{i+}\otimes\mu_t^{i+})
\Big(\big\{ (x,y)\,;~~|x-y|\le \delta\big\}\Big)\,dt
=\O(1)\cdot \Big(\ln(2+\tau)+|\ln\delta|\Big)\delta\cdot\tv\{\bar u\}\,.
\eqno(2.25)$$
\v
\n{\bf Proof.}
By Proposition 1 and the previous comparison lemmas,
for every $i=1,\ldots,n$ the integral on the left hand side of (2.25)
has the same order of magnitude as in the case of a 
solution with a single
centered rarefaction wave, of magnitude $\sigma\doteq\tv\{\bar u\}<1$.
Looking back at Example 1, from (2.3) we thus obtain
$$\eqalign{\int_0^\tau (\mu_t^{i+}\otimes\mu_t^{i+})
\Big(\big\{ (x,y)\,;~~|x-y| \le 2\delta\big\}\Big)
&=\O(1)\cdot 
\delta\sigma\left(1+
\ln{\sigma\tau\over 2\delta}\right)\cr
&=\O(1)\cdot
\Big(\ln(2+\tau)+|\ln\delta|\Big)\delta\cdot\tv\{\bar u\}\,.\cr}\eqno(2.26)$$
\endproof
\v
\n{\bf Remark 3.} 
All of the above estimates refer to an exact solution $u$ of (1.1).
If $u_\nu\to u$ is a convergent sequence of front tracking approximations,
the corresponding measures
of $i$-waves in $u_\nu(t,\cdot)$ converge weakly:
$\mu^i_{\nu,t}\wto \mu^i_t$ for all $i=1,\ldots,n$ and $t\geq 0$.
Unfortunately, this does not guarantee the weak convergence of the
signed measures
$$\mu^{i+}_{\nu,t}\wto \mu^{i+}_t\,,\qquad\qquad 
\mu^{i-}_{\nu,t}\wto \mu^{i-}_t\,.\eqno(2.27)$$
For example (fig.~7), on a fixed interval $[a,b]$
every $u_\nu$ might contain an alternating sequence
of small positive and negative waves, that
cancel only in the limit as $\nu\to\infty$.
However, by a small modification of these 
front tracking
solutions $u_\nu$ one can achieve the weak convergence
(2.27) for each $t$ in a discrete set of times 
$\{j\tau/N\,;~j=0,1,\ldots,N\}$,
with $N>\!>\ve^{-1}$.
As a result, we obtain an arbitrarily accurate 
front tracking approximation (still called $u$)
satisfying an estimate entirely analogous to (2.25), namely
$$\sum_{i=1}^n\int_0^\tau \bigg(\sum_{\alpha,\beta\in\Rar_i\,,~
|x_\alpha-x_\beta|\le 8\sqrt\ve}|\sigma_\alpha\sigma_\beta|\bigg)\,dt
=\O(1)\cdot \Big(\ln(2+\tau)+|\ln\ve|\Big)\sqrt\ve\cdot\tv\{\bar u\}\,,
\eqno(2.28)$$
where we replace the $\delta $ in (2.25) by $8\sqrt\ve$ for the
application in Section 4.
Here $\Rar_i$ denotes the set of
rarefaction fronts of the $i$-th family and summation is
over all possible pairs, including the case where the two indices 
$\alpha,\beta$ coincide.

\vsk
\n{\medbf 3 - Estimates on shock fronts} 
\v
We begin by estimating the sum in (1.9).
The approximation $v$ is discontinuous precisely at those
times $t_i$ where an interaction occurs involving a large shock.
Indeed, at such times the left and right states $u_\alpha^-$, $u_\alpha^+$
across a large shock located at $x=x_\alpha$ suddendly change. 
As a consequence, 
the viscous shock profile connecting these two states is modified.
The two smooth functions $v(t_i-)$ and $v(t_i)$ will thus be different
over the interval $[x_\alpha-\sqrt\ve,\, x_\alpha+\sqrt\ve]$.
To estimate the $\L^1$ norm of this difference, 
the following elementary observation is
useful.  Given a smooth function $\phi=\phi(\sigma,\sigma')$, its size
satisfies the bounds:
$$\hbox{if}~~~\phi(\sigma,0)=0~~~\hbox{for all}~\sigma,~~\hbox{then}~~
\phi(\sigma,\sigma')=\O(1)\cdot |\sigma'|\,,$$
$$\hbox{if}~~~\phi(\sigma,0)=\phi(0,\sigma')=0~~~\hbox{for 
all}~\sigma,\sigma',~~\hbox{then}~~
\phi(\sigma,\sigma')=\O(1)\cdot |\sigma\,\sigma'|\,.$$

We now distinguish various cases.
\v
\n{\bf 1.} At time $t_i$ a new large shock is created, say of strength
$|\sigma_\alpha|\geq \rho/2$.
In this case, since the new viscous shock profile is
inserted on an interval of length $2\sqrt\ve$, we have
$$\big\|v(t_i)-v(t_i-)\big\|_{\L^1} 
=\O(1)\cdot\sqrt\ve\,|\sigma_\alpha|\,.$$
According to our construction, every large shock not present at time
$t=0$ must grow from a strength $<\rho/2$ up to a strength
$\geq \rho$ at some later time $\tau$.
Therefore, the sum of the strengths of all
large shocks, at the time when $t_i$ when
they are created, is $\O(1)\cdot \delta_0$, where
$\delta_0\doteq\tv\{\bar u\}$.
The total contribution due to these terms is thus 
$\O(1)\cdot\sqrt\ve\,\delta_0$.
\v
\n{\bf 2.} At time $t_i$ a large shock is terminated.
Since every large shock must have
strength $\geq\rho$ at some time and is terminated when
its strength becomes $<\rho/2$, every such case
involves an amount of interaction and cancellation $\geq \rho/2$.
Therefore, the total contribution of these terms to the sum in (1.9)
is again $\O(1)\cdot\sqrt\ve\,\delta_0$.
\v
\n{\bf 3.}
A front $\sigma_\beta$ of a different family crosses one large shock
$\sigma_\alpha$.
In this case we have
$$\big\|v(t_i)-v(t_i-)\big\|_{\L^1} =\O(1)\cdot \sqrt\ve\,|\sigma_\alpha|
\,|\sigma_\beta|\,.$$
These terms are thus controlled by the decrease in the interaction
potential $Q(u)$.  Their total sum is 
$\O(1)\cdot \sqrt\ve\,\delta_0^2$.
\v
\n{\bf 4.}
A small front $\sigma_\beta$ of the same family impinges on the large shock
$\sigma_\alpha$.
In this case we have
$$\big\|v(t_i)-v(t_i-)\big\|_{\L^1} =\O(1)\cdot \sqrt\ve\,|\sigma_\beta|
\,,$$
Since any small front can join at most one large shock of the same
family, the total contribution of these terms
is $\O(1)\cdot \sqrt\ve\,\delta_0$.
\v
\n{\bf 5.}
Two large $k$-shocks of the same family, say of strengths $\sigma_\alpha,\,
\sigma_\beta$, merge together.
In this case
$$\big\|v(t_i)-v(t_i-)\big\|_{\L^1} =\O(1)\cdot \sqrt\ve\,\min\big\{
|\sigma_\alpha|,
\,|\sigma_\beta|\big\}
\,.$$
As will be shown in (3.23), all these interactions are controlled 
by the decrease in
a suitable functional $Q^\sharp(u)$ by noticing that
$|\sigma_\alpha|, |\sigma_\beta| >2\sqrt\ve|\ln\ve|$.
The sum of all these terms is thus found to be 
$\O(1)\cdot \delta_0\sqrt\ve\,|\ln\ve|\,$.
\v
\n Putting together all these five cases, one obtains the bound (1.9).
\vs
Next, we need to estimate the running error in (1.8)
related to the big shocks, namely
$$E_\BS\doteq\int_0^\tau\sum_{\alpha\in \BS(t)}\int_{x_\alpha-\sqrt\ve}
^{x_\alpha+\sqrt\ve} \big|v_t+A(v)v_x-\ve v_{xx}\big|\,dxdt\,.\eqno(3.1)$$
Here the summation ranges over all big shocks in $v(t,\cdot)$.

We first consider the simplest case, where the interval
$$I_\alpha(t)\doteq \big[x_\alpha(t)-2\sqrt\ve\,,~x_\alpha(t)+2\sqrt\ve
\big]\eqno(3.2)$$
does not contain any other wave-front.
In this case, observing that
$$\big[ A(\omega_\alpha^\ve(s))-\lambda_\alpha\big] 
{\partial\over\partial s}\omega_\alpha^\ve(s)-\ve\,
{\partial^2\over\partial s^2}\omega_\alpha^\ve(s)=0\,,$$
and recalling (1.19)-(1.20), the error relative to the shock
at $x_\alpha$ can be written as
$$\eqalign{E_\alpha(t)&=\left(\int_{-\sqrt\ve}^{-\sqrt\ve/2}
+\int_{\sqrt\ve/2}^{\sqrt\ve}\right)
\bigg\{
\Big[ A\big(\omega_\alpha^\ve(\phi(\xi))\big)-\lambda_\alpha\Big] 
{\partial\over\partial s}\omega_\alpha^\ve(\phi(\xi))\,\phi'(\xi)\cr
&\qquad\qquad-\ve\,
{\partial\over\partial s}\omega_\alpha^\ve(\phi(\xi))\phi''(\xi)
-\ve\,
{\partial^2\over\partial s^2}\omega_\alpha^\ve(\phi(\xi))\big(\phi'(\xi)
\big)^2\bigg\}\,d\xi\,.\cr}\eqno(3.3)$$
Using the bounds
$$\left|{\partial\over\partial s}\omega_\alpha^\ve(s)\right|=\O(1)
\cdot {|\sigma_\alpha|^2\over\ve}\,e^{-|s\,\sigma_\alpha|/\ve}\,,\qquad
\left|{\partial^2\over\partial s^2}\omega_\alpha^\ve(s)\right|\,=\O(1)
\cdot {|\sigma_\alpha|^3\over\ve^2}\,e^{-|s\,\sigma_\alpha|/\ve},
\eqno(3.4)$$
{}from (1.20) we deduce
$$
\eqalign{E_\alpha(t)&=\O(1)\cdot\int_{\sqrt\ve/2}^{\sqrt\ve}
\exp\left\{ -{|\sigma_\alpha|\over\ve}\,\phi(\xi)\right\}\cdot
\left({|\sigma_\alpha|^2\over\ve}\,\phi'(\xi)+
{|\sigma_\alpha|^3\over\ve}\,\phi''(\xi)\right)\,d\xi\cr
&=\O(1)\cdot\int_{\sqrt\ve/2}^{\sqrt\ve}
\exp\left\{ -{|\sigma_\alpha|\over 4(\sqrt\ve-\xi)}\right\}\, 
{|\sigma_\alpha|^3\over (\sqrt\ve-\xi)^3}\,d\xi\cr
&=\O(1)\cdot\int_{2/\sqrt\ve}^\infty 
\exp\left\{ -{|\sigma_\alpha|\,s\over 4}
\right\}\, |\sigma_\alpha|^3\,s^3\,{ds\over s^2}\cr
&=\O(1)\cdot |\sigma_\alpha|\,
\exp\left\{ -{|\sigma_\alpha|\over 2\sqrt\ve}\right\}\,
\left( 1+{2|\sigma_\alpha|\over\sqrt\ve}\right).
\cr}
$$
Since by assumption $|\sigma_\alpha|\geq \rho/2=2\sqrt\ve\,|\ln\ve|$, 
the above 
estimate implies
$$E_\alpha(t)=\O(1)\cdot\ve
\big( 1+|\ln\ve|\big)\,|\sigma_\alpha|\,.
\eqno(3.5)$$
In the general case, our error estimate
must also take into account the presence of other wave-fronts
within the intervals $I_\alpha(t)$.  Indeed, for every point $x_\alpha$
where large shock is located, we have
$$
\eqalign{
E_\alpha(t)&\doteq
\int_{x_\alpha-\sqrt\ve}
^{x_\alpha+\sqrt\ve} \big|v_t+A(v)v_x-\ve v_{xx}\big|\,dxdt\cr
&=
\O(1)\cdot\ve
\big( 1+|\ln\ve|\big)\,|\sigma_\alpha|
+
\O(1) \left(\sum_{x_\beta, x_\gamma\in I_\alpha(t),
|x_\beta-x_\gamma|\le 2\sqrt\ve}
|\sigma_\beta\sigma_\gamma|-\sum_{x_\theta \in I_\alpha(t),
\theta \in\BS}|\sigma_\theta|^2\right).
\cr}
$$
\v
In the following, we introduce three different functionals, 
which account for:

\i{$\bullet$} products $|\sigma_\alpha\sigma_\beta|$ 
of fronts of different families,

\i{$\bullet$} products $|\sigma_\alpha\sigma_\beta|$
where $\sigma_\alpha$ is a large shock and 
$\sigma_\beta$ is a rarefaction of the same family, 

\i{$\bullet$} products $|\sigma_\alpha\sigma_\beta|$ 
of shocks the same family.
\v
By combining these three, we form a functional $\Hat Q(u)$ such that
the map $t\mapsto \Hat Q\big(u(t)\big)$ is non-increasing
except at times where a new large shock is introduced. 
Moreover, the total increase in this functional
at times where large shocks are created will be shown to be
$\O(1)\cdot 
\sqrt\ve|\ln\ve|\tv\{\bar{u}\}$. 
\v
We begin by defining
$$Q^\flat(u)\doteq 
\sum_{k_\beta\not= k_\alpha} 
W^\flat_{\alpha\beta}|\sigma_\alpha \sigma_\beta|
\,.\eqno(3.7)$$
where the sum extends over all couples of fronts of different families
(small shocks, big shocks, rarefactions).
The 
weights $W^\flat_{\alpha\beta}\in [0,1]$ 
are defined as follows. If $k_\beta<k_\alpha$, 
then
$$W^\flat_{\alpha\beta}\doteq \left\{
\eqalign{0\qquad &\qquad\hbox{if}\qquad x_\beta<x_\alpha-2\sqrt\ve\,,\cr
\qquad {1\over 2}+{x_\beta-x_\alpha\over 4\sqrt\ve}
&\qquad\hbox{if}\qquad x_\beta\in [x_\alpha-2\sqrt\ve\,,~x_\alpha+2\sqrt\ve
]\,,\cr
1\qquad &\qquad\hbox{if}\qquad x_\beta>x_\alpha+2\sqrt\ve\,.\cr}
\right.$$
If instead $k_\beta>k_\alpha$, we set
$$W^\flat_{\alpha\beta}\doteq \left\{
\eqalign{1\qquad &\qquad\hbox{if}\qquad x_\beta<x_\alpha-2\sqrt\ve\,,\cr
\qquad {1\over 2}-{x_\beta-x_\alpha\over 4\sqrt\ve}
&\qquad\hbox{if}\qquad x_\beta\in [x_\alpha-2\sqrt\ve\,,~x_\alpha+2\sqrt\ve
]\,,\cr
0\qquad &\qquad\hbox{if}\qquad x_\beta>x_\alpha+2\sqrt\ve\,.\cr}\right.$$
By strict hyperbolicity, we expect that the functional $Q^\flat$ 
will be decreasing in time.  Indeed, its rate of decrease
dominates the sum 
$$\sum_{k_\alpha\not= k_\beta,\,|x_\alpha-x_\beta|<2\sqrt\ve}
|\sigma_\alpha\sigma_\beta|\,,$$
containing products of nearby waves of different families.
\v

\midinsert
\vskip 10pt
\centerline{\hbox{\psfig{figure=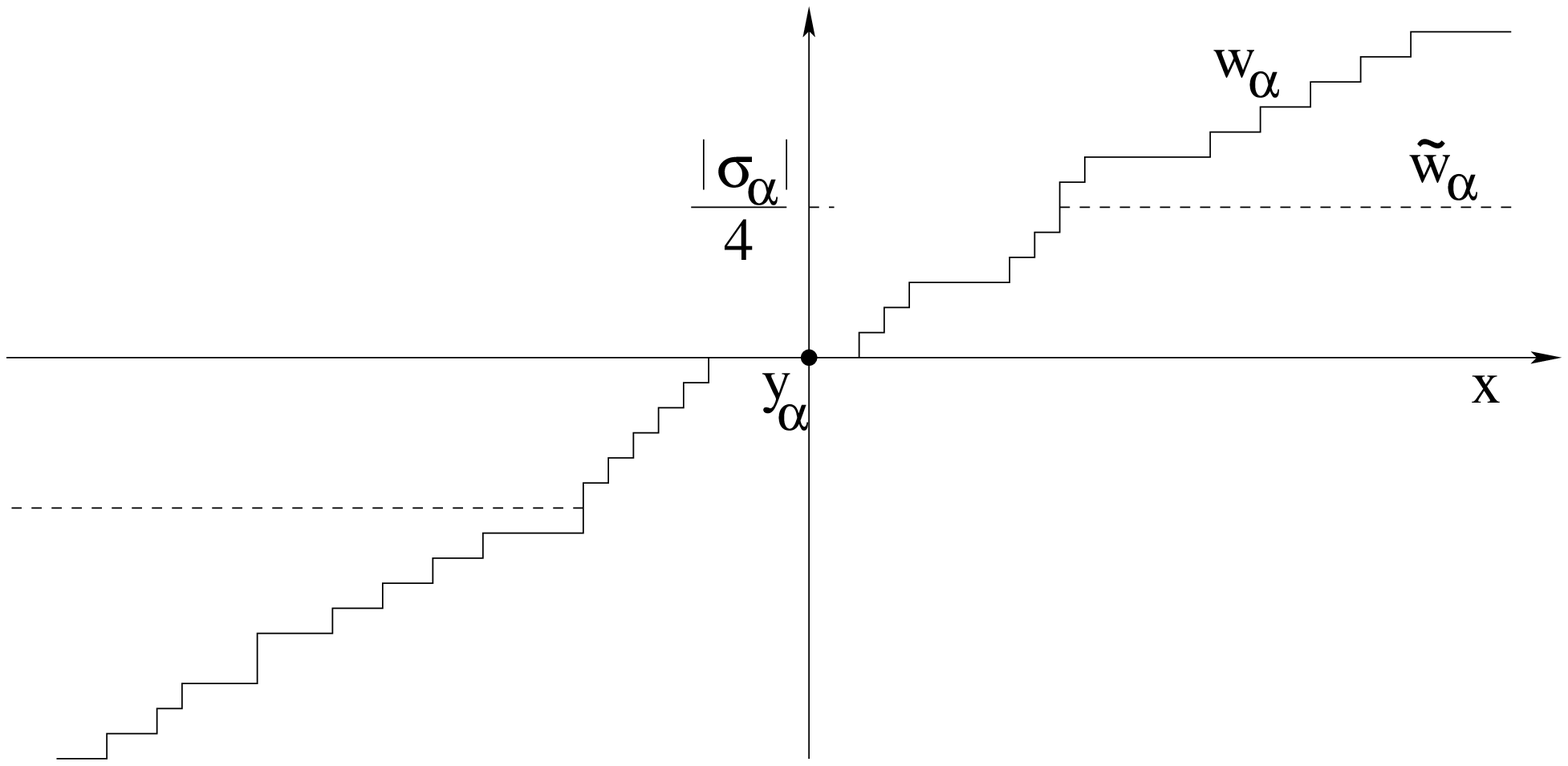,width=10cm}}}
\centerline{\hbox{figure 4}}
\vskip 10pt
\endinsert

Next, given a big shock  
$\sigma_\alpha$ of the $k_\alpha$-th family located at $x_\alpha$, 
we write:
\v
$\Rar_\alpha$ to denote the set of all rarefaction fronts 
of the same family $k_\alpha$,

$\S_\alpha$ to denote the set of all shock fronts 
of the same family $k_\alpha$.
\v
To control the interaction between large shocks and rarefactions
of the same family, we define the weight
$$W^\nat_{\alpha}(x)\doteq \min\left\{ {1\over 2}+{|x -x_\alpha|\over 
4\sqrt\ve}\,,~ 1\right\}\eqno(3.8)$$
and the function (fig.~4)
$$w_\alpha(x)\doteq \left\{\eqalign{
&\sum_{\beta\in\Rar_\alpha,\,
x_\beta\in [x, \,x_\alpha]} 
(-\sigma_\beta)\qquad\quad \hbox{if}\quad 
x<x_\alpha\,,\cr
&\cr
&\sum_{\beta\in\Rar_\alpha,\,
x_\beta\in [x_\alpha,\,x]} 
~~\sigma_\beta~~\qquad\quad \hbox{if}\quad 
x>x_\alpha\,.\cr}\right.$$
Calling 
$$\tilde w_\alpha(x)\doteq \cases{ -|\sigma_\alpha|/4\qquad &if
\quad $w_\alpha(x)<-|\sigma_\alpha|/4\,$,\cr
&\cr
~w_\alpha(x)\qquad &if
\quad $\big |w_\alpha(x)\big|\leq |\sigma_\alpha|/4\,$,\cr
&\cr
|\sigma_\alpha|/4\qquad &if
\quad $w_\alpha(x)>|\sigma_\alpha|/4\,$,\cr}$$
we then define
$$Q^\nat(u)\doteq\sum_{\alpha\in\BS}
\int W^\nat_\alpha(x)\, D_x \tilde w_\alpha\,.\eqno(3.9)$$
By using the function with cut-off $\tilde w_\alpha$,
instead of $w_\alpha$, in (3.9) we are taking into account only
the rarefaction fronts $\sigma_\beta$
of the same family $k_\alpha$, such that the total amount of
rarefactions inside the interval 
$[x_\alpha, x_\beta]$ is $\leq |\sigma_\alpha|/4$.
If no other fronts of different families are present,
this guarantees that all these rarefactions $\sigma_\beta$
are strictly
approaching the big shock $\sigma_\alpha$. Indeed,
the difference in speed is $|\dot x_\beta-\dot x_\alpha|
\geq |\sigma_\alpha|/4$. 
As a result, the functional $Q^\nat(u) $ will be strictly decreasing.
On the other hand, if the interval $[x_\alpha, x_\beta]$
also contains
waves of different families, the above estimate may fail.
In this case, however, the decrease in the functional $Q^\flat(u)$
 compensates the possible increase in $Q^\nat(u)$.

\midinsert
\vskip 10pt
\centerline{\hbox{\psfig{figure=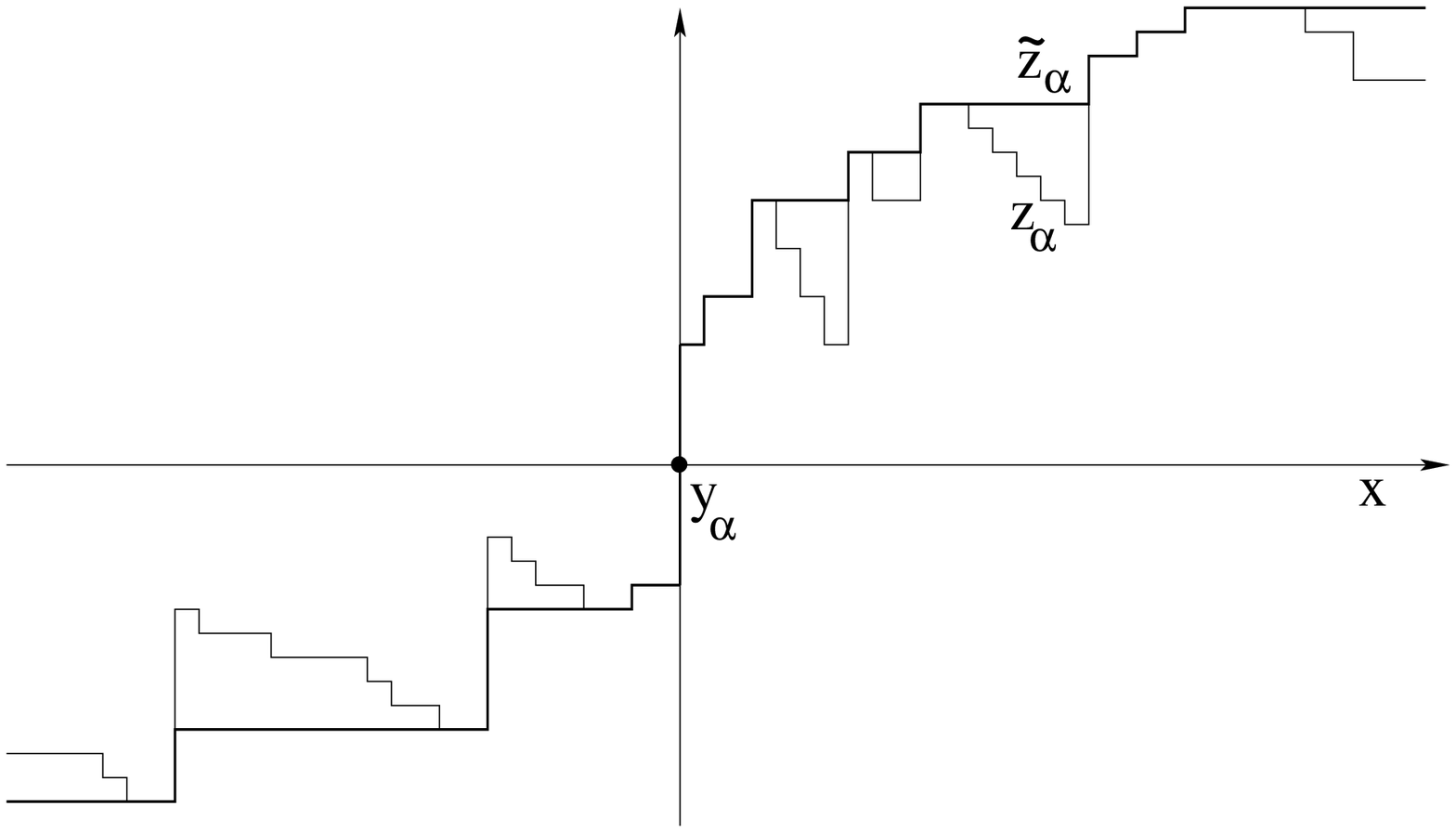,width=10cm}}}
\centerline{\hbox{figure 5}}
\vskip 10pt
\endinsert

Finally, to control the interactions among shocks
of the same family, 
for each shock front $\sigma_\alpha$ (of any size, big or small)
located at $x_\alpha$,
we begin by defining (fig.~5)
$$z_\alpha(x)\doteq\left\{
\eqalign{-{|\sigma_\alpha|\over 2}-
\sum_{\beta\in\S_\alpha,~x<x_\beta
< x_\alpha} |\sigma_\beta|+\sum_{\beta\in\Rar_\alpha,~x<x_\beta
< x_\alpha} 3\sigma_\beta
~~\qquad &\hbox{if} \quad x< x_\alpha\,,\cr
{|\sigma_\alpha|\over 2}
+\sum_{\beta\in\S_\alpha,~x_\alpha<x_\beta
<x} |\sigma_\beta|-\sum_{\beta\in\Rar_\alpha,~x_\alpha<x_\beta<x
} 3\sigma_\beta
~~\qquad &\hbox{if} \quad x> x_\alpha\,.\cr
}\right.$$
Then we set
$$\tilde z_\alpha(x)=\cases{\min\big\{ z_\alpha(x')\,;~~x<x'<x_\alpha\big\}
\qquad &if\qquad $x<x_\alpha\,$,\cr
&\cr
\max\big\{ z_\alpha(x')\,;~~x_\alpha<x'<x\big\}
\qquad &if\qquad $x>x_\alpha\,$.\cr}$$
Notice that $\tilde z_\alpha$ is a non-decreasing, piecewise constant 
function, with $(x-x_\alpha)\,\tilde z_\alpha(x)>0$
for $x\not= x_\alpha$.

Using the weights
$$W_\alpha^\sharp(x)\doteq \cases{\big[\ve-\tilde z_\alpha(x-)\big]^{-1}
\qquad &if\quad $x<x_\alpha\,$,\cr
\big[\ve +\tilde z_\alpha(x+)\big]^{-1}
\qquad &if\quad $x>x_\alpha\,$,\cr}$$
we now define
$$
Q^\sharp(u)
\doteq\sum_{\alpha\in\S} |\sigma_\alpha|\,
\int W_\alpha^\nat(x)\,W_\alpha^\sharp(x) D_x\tilde z_\alpha\,.
\eqno(3.10)$$
Notice that in this case the summation runs over all shock fronts.
If $\sigma_\beta$ is a shock
located at $x_\beta$, then 
$\big[W_\alpha^\sharp(x_\beta)\big]^{-1}$
roughly 
describes the amount of shock waves inside the interval $[x_\alpha, x_\beta]$
in excess of three times the amount of rarefactions.
If the interval $[x_\alpha, x_\beta]$
does not contain waves of other families
and the function $x\mapsto \tilde z_\alpha(x)$ 
has a jump at at $x=x_\beta$, then the two shocks $\sigma_\alpha$,
$\sigma_\beta$ are strictly approaching, hence the functional
$Q^\sharp(u)$ will decrease.  
On the other hand, if 
waves of different families are present, the above estimate may fail.
In this case, however, the decrease in the functional $Q^\flat(u)$
 compensates the possible increase in $Q^\sharp(u)$.

In the definition of $z_\alpha$, notice that
the strength of rarefactions is 
multiplied by 3, to make sure that couples of shocks
$\sigma_\alpha$, $\sigma_\beta$ entering the definition of
$Q^\sharp(u)$ are always approaching each other 
(except for the presence of fronts of different families in 
between).
An example is shown in
fig.~6, where  two nearby shocks move apart from each other
because there are sufficiently many rarefaction waves in the middle.
Because of the factor $3$, the function $x\mapsto \tilde z_\alpha(x)$ 
will be constant at the point $x_\beta$. 
Hence the product $|\sigma_\alpha\sigma_\beta|$
will not appear within the definition of $Q^\sharp(u)$.

\midinsert
\vskip 10pt
\centerline{\hbox{\psfig{figure=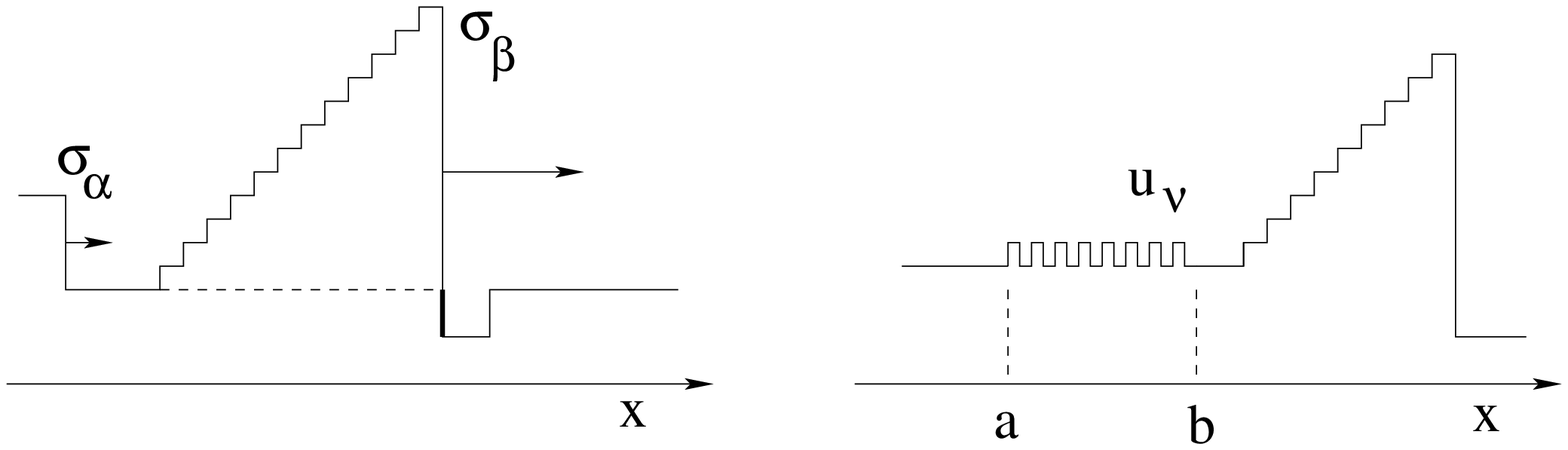,width=14cm}}}
\centerline{\hbox{figure 6~~~~~~~~~~~~~~~~~~~~~~~~~~~~
~~~~~~~~~~~~~~~~~~~~~~figure 7}}
\vskip 10pt
\endinsert

\v
We now consider the composite functional
$$\Hat Q(u)\doteq \sqrt\ve\,|\ln\ve|\cdot 
\Big(C_1\,\Ups(u)+C_2\,Q^\flat(u)+ C_3Q^\nat(u)\Big)+\sqrt\ve \,Q^\sharp(u)\,.
\eqno(3.11)$$
Here
$$\Ups(u)\doteq V(u)+C_0\,Q(u)\eqno(3.12)$$
is a quantity which is decreasing at every interaction time. 
Its decrease dominates both
the amount of interaction and of cancellation 
in the front tracking solution $u$.
Observe that
$$\Hat Q(u)=\O(1)\cdot \sqrt\ve|\ln\ve|\,\tv\{u\}\,.\eqno(3.13)$$
Indeed, by the definition of $W^\sharp_\alpha(x)$,  we have
$$
\int W_\alpha^\sharp(x) D_x\tilde z_\alpha=\O(1)\cdot
\int_{0}^{\tv\{u\}} ~
{1\over {s+\ve}}\,ds 
=\O(1)\cdot |\ln\ve|\,.\eqno(3.14)$$
Using (3.14), it is now clear that
$$Q^\sharp(u)=\O(1)\cdot |\ln\ve|\,\tv\{ u\}\,,\qquad
\quad\Ups(u),~Q^\flat(u),~
Q^\nat(u)
=\O(1)\cdot\tv\{u\}\,.\eqno(3.15)$$
The bound on (3.11) now follows from (3.15).
\v
\n{\bf Lemma 5.} 
{\it For a suitable choice of the constants $C_1>\!>C_2>\!>C_3>\!> 1$, if
$\tv\{u\}$ remains small, then
at each time $t^*$ where an
interaction occurs the following holds.
If a new large shock of strength $|\sigma_\alpha|> 2\sqrt\ve|\ln\ve|$
is created, then
$$\Delta\Hat Q\doteq \hat Q(\tau+)-\hat Q(\tau-)=
\O(1)\cdot \sqrt\ve|\ln\ve||\sigma_\alpha|\,.\eqno(3.16)$$
If no large shock is created, then
$$\Delta\Hat Q\leq 0\,.\eqno(3.17)$$
}
\v
\n{\bf Proof.} 
Notice that the weight $W^\flat_{\alpha,\beta}$ is 
always $\leq 1$. For a newly created large shock 
$\sigma_\alpha$, the increase in
the functional $Q^\flat(u)$ can be estimated as
$$
\Delta Q^\flat(u)=\O(1)\cdot |\sigma_\alpha|\,\tv\{u\}.\eqno(3.18)
$$
Similarly, since $W^\nat_\alpha\leq 1$, 
it is clear that the increase of 
$Q^\nat(u)$ due to a new large shock $\sigma_\alpha$ is
$$
\Delta Q^\nat(u)=\O(1)\cdot |\sigma_\alpha|.\eqno(3.19)
$$
The estimate on the increase of the  functional $Q^\sharp(u)$ is
different. In this case, the integral
 $$\int W^\nat_\alpha(x)W^\sharp_\alpha(x)
D_x\tilde{z}_\alpha$$ is bounded by 
$$
\O(1)\cdot\int_0^{\tv\{u\}} {1\over{\ve+ x}}\,dx=\O(1)|\ln\ve|.
$$
 Hence,
$$\Delta Q^\sharp(u) =\O(1)\cdot |\sigma_\alpha|
\,|\ln\ve|.\eqno(3.20)$$
Together, (3.18)-(3.20) imply (3.16).

Next, we prove (3.17). 
Assume that at time $t^*$ an interaction occurs without 
the introduction of any new
large shock. We will show that  the functional
$\Hat{Q}(u(t))$ decreases. 

First we look at the change in $Q^\flat(u)$ and $Q^\nat(u)$. Since
the weights $W^\flat$ and $W^\nat$ are uniformly bounded, it is
straightforward to check that the change in these two functionals
at time $t^*$ is bounded by a constant times the decrease in the
Glimm functional $\Ups(u(t))$ in (3.12).  Hence, by choosing
$C_1>\!> C_2>\!>C_3$, the quantity
$$
C_1\Ups(u) +C_2 Q^\flat(u) +C_3Q^\nat(u)
$$
is not increasing in time.

The analysis of  $Q^\sharp(u)$ is a bit harder. We will
show that the change of $Q^\sharp(u)$ at the interaction time $t^*$ 
is of the same order of magnitude as $|\ln\ve|\big|\Delta \Ups(u)\big|$.
Here and in the following, $\Delta \Ups$ denotes the change in
$\Ups(u(t))$ across the interaction time. 
As a preliminary,
we notice a 
basic property of the weight function $W^\sharp_\alpha(x)$.
 For any fixed location
$x=x_0$, we have
$$
\sum_{\alpha\in S} |\sigma_\alpha|\,W^\sharp_{\alpha}(x_0)
=\O(1)\cdot |\ln\ve|,\eqno(3.21)
$$
$$
\eqalign{
\sum_{\alpha\in S,\, x(\alpha)<x_0} |\sigma_\alpha|
\int_{x_0}^\infty  (W^\sharp_\alpha(x))^2 D_x\tilde z_\alpha
&+
\sum_{\alpha\in S,\, x(\alpha)>x_0} |\sigma_\alpha|
\int^{x_0}_{-\infty} \big(W^\sharp_\alpha(x)\big)^2 D_x\tilde z_\alpha
\cr
&=\O(1)\cdot|\ln\ve|.\cr}\eqno(3.22)
$$
The proof of the estimates in (3.21) and (3.22) is
straightforward by noticing that the functions $f(x)={{1}\over {x+\ve}}$,
$g(x)={{1  }\over{(x+\ve)^2   } } $
are convex and bounded away from zero for $x\ge 0$.
And the left hand sides of (3.21) and (3.22)
are   bounded by the following single and  double
integrals respectively:
$$
\int_0^{\tv\{u\}} f(x)=\O(1)|\ln\ve|,\qquad\qquad
\int_0^{\tv\{u\}}\int_x^{\tv\{u\}} g(y)dydx=\O(1)\cdot |\ln\ve|.
$$

Now we are ready to estimate the change in $Q^\sharp(u(t))$ at time
$t^*$. Note that, in some cases, it is possible
that the interaction does not change the functional. In the following, we
will consider the case where 
$Q^\sharp(u)$ does change across the interaction.
Depending on the types and families
of the waves involved in the interaction, we have the following four
cases.

\v

\n{\bf 1.} Two shocks of the same family interact.
Let $\beta_1$ and $\beta_2$ be the two interacting shocks, say of the
$i$-th family, and call 
$\beta$ the outgoing  $i$-shock. We also let
$\alpha_1$ and $\alpha_2$ be any two shock fronts on the left and right
of the interaction point respectively, so that
$x_{\alpha_1}<x_\beta<x_{\alpha_2}$ at time $t=t^*$.
For any shock front $\sigma_\alpha$ at time $t^*$, set
$$
Q^\sharp_{\alpha}=
|\sigma_{\alpha}|\int W^\nat_{\alpha}(x)
W^\sharp_{\alpha}(x) \,D_x\tilde z_{\alpha}.
$$
Observe that
$$
Q^\sharp_{\beta} -\big(Q^\sharp_{\beta_1} +Q^\sharp_{\beta_2}\big)
\leq - {|\sigma_{\beta_1}\sigma_{\beta_2}|\over {|\sigma_{\beta_1}
+|\sigma_{\beta_2}| +\ve}}+
\O(1)\cdot|\ln\ve|\,|\Delta \Ups|\,.\eqno(3.23)
$$
Indeed,
$$
\sigma_\beta=\sigma_{\beta_1}+\sigma_{\beta_2} + \O(1)\cdot
|\Delta \Ups|.
$$
Moreover, recalling (3.8), we see that after the interaction
we lose the term
$$W^\nat_{\beta_1}(x_{\beta_2})(
W^\sharp_{\beta_1}(x_{\beta_2}) +
W^\sharp_{\beta_2}(x_{\beta_1})) |\sigma_{\beta_1}\sigma_{\beta_2}|
\geq {|\sigma_{\beta_1}\sigma_{\beta_2}|\over 
|\sigma_{\beta_1}|+|\sigma_{\beta_2}|+\ve }.$$

Notice that $W^\sharp_{\alpha_1}(x)$ (respectively $W^\sharp_{\alpha_2}(x)$)
does not change across the interaction for $ x<x_\beta$
($x>x_\beta$).
The change in the $Q^\sharp_{\alpha_i}$, $i=1,\, 2$, 
can be estimated as follows. 
When $\alpha_i$ is of the same family of $\beta_j$, 
$(i,j=1,2)$, by (3.22) we have
$$
\eqalign{
\Delta Q^\sharp_{\alpha_1}
&=
\bigg({{ |\sigma_{\alpha_1}\sigma_\beta| }\over{  \ve +|\sigma_{\alpha_1}|
+|\sigma_\beta|+I } }
 -{ {|\sigma_{\alpha_1}\sigma_{\beta_1}| }\over{\ve
 +|\sigma_{\alpha_1}|+|\sigma_{\beta_1}|
+I   } } -{{|\sigma_{\alpha_1}\sigma_{\beta_2}|  }
\over{\ve +|\sigma_{\alpha_1}|+|\sigma_{\beta_1}|
+|\sigma_{\beta_2}|+I    }}\bigg) W^\nat_\alpha(x(\beta))\cr
&\qquad\qquad
+\O(1)\cdot |\Delta \Ups||\sigma_{\alpha_1}|\int_{x_\beta}^\infty
(W^\sharp_{\alpha_1})^2(x) D_x\tilde z_{\alpha_1}\cr
&\le\O(1)\cdot |\sigma_{\alpha_1}||\Delta \Ups| W_{\alpha_1}^\sharp(x_\beta)+
\O(1)\cdot|\Delta \Ups||\sigma_{\alpha_1}|\int_{x_\beta}^\infty
(W^\sharp_{\alpha_1})^2(x) D_x\tilde z_{\alpha_1}\,.
\cr}
$$
Here and in the following, we assume that the whole
strength of
$\beta_i$, $i=1,\, 2$ and of $\beta$ appear in the functional
$Q^\sharp_{\alpha_1}$. Moreover, $I$ 
represents the sum of the strengths of the $i$-shocks
between $\beta_1$ and $\alpha_1$ that appear in 
$Q^\sharp_{\alpha_1}$. The other cases when part or none of the
above wave stengths appears in $Q^\sharp_{\alpha_1}$ can be treated
similarly. 

By summing over $\alpha_1$ and using (3.21) and (3.22), we find that
the total change  of $ Q^\sharp_{\alpha_1}$ is $\O(1)\cdot |\ln\ve|
|\Delta \Ups|$.  A similar estimate holds for $\alpha_2$.

Now  consider two shock fronts
$\alpha_1$, $\alpha_2$ of the $j$-th
family, with $j\neq i$. Notice that the change of the weight
function $W^\sharp_{\alpha_i}(x)$, $i=1,\, 2$, is at most of the order
of $(W^\sharp_{\alpha_i}(x))^2|\Delta \Ups|$
when $x$ lies on the opposite side of $x_{\alpha_i}$ w.r.t.
 $x_\beta$.  Together with (3.22) this yields 
$$
\eqalign{
\sum_{\alpha_i,\, i=1,2}|\Delta Q^\sharp_{\alpha_i}|
&=\O(1)\cdot\left(
\sum_{\alpha_1}|\Delta \Ups||\sigma_{\alpha_1}|\int_{x_\beta}^\infty
(W^\sharp_{\alpha_1})^2(x) D_x\tilde z_{\alpha_1}
+\sum_{\alpha_2}|\Delta \Ups||\sigma_{\alpha_2}|\int^{x_\beta}_{-\infty}
(W^\sharp_{\alpha_2})^2(x) D_x\tilde z_{\alpha_2}
\right)\cr
&=\O(1)\cdot|\ln\ve||\Delta \Ups|\,.\cr
}
$$

If $\gamma$ is a newly created shock of the $j$-th family, then
the new term $Q^\sharp_{\gamma}$ has size 
$\O(1)\cdot|\sigma_\gamma||\ln\ve|$. 
Hence, the total sum of these new terms over
$\gamma$ is of $\O(1)|\ln\ve ||\Delta \Ups|$. And this completes the
discussion on this case.

\v

\n{\bf 2.} Interaction of a shock with a rarefaction front of the same 
family.  Let
$\beta_1$ and $\beta_2$ be a shock and a rarefaction front of the 
$i$-th family,  
interacting at time $t^*$.  

First, consider the case where the shock
$\beta_1$ is completely cancelled and hence the decrease in 
$\Ups(u)$ is of the same order as $\beta_1$. 
In this case the term $Q^\sharp_{\beta_1}$ disappears after the 
interaction.
Let $\alpha_1$ and $\alpha_2$ be shock waves of the $j$-th
family on the left and right of the location of interaction.
For both cases when $i=j$ and $i\ne j$, by (3.22) we have
$$\eqalign{
\sum_{\alpha_i,\, i=1,2}
|\Delta Q^\sharp_{\alpha_i}|&=\O(1)\cdot \left(
\sum_{\alpha_1}|\Delta \Ups||\sigma_{\alpha_1}|\int_{x_\beta}^\infty
(W^\sharp_{\alpha_1})^2(x) D_x\tilde z_{\alpha_1}
+\sum_{\alpha_2}|\Delta \Ups||\sigma_{\alpha_2}|\int^{x_\beta}_{-\infty}
(W^\sharp_{\alpha_2})^2(x) D_x\tilde z_{\alpha_2}
\right)\cr
&=\O(1)\cdot |\ln\ve||\Delta \Ups|.\cr
}
$$ 
The same argument applies to the change in $ Q^\sharp_\gamma$,
related to
the newly created shock $\gamma$ of the
$j$-th family, when $j\neq i$. In this case, the total
change in $Q^\sharp$ is again $\O(1)\cdot|\ln\ve||\Delta \Ups|$.

In the case where the interaction produces an outgoing $i$-shock
$\bar{\beta}_1$, so that
the rarefaction $\beta_2$ is completely cancelled,
the analysis is as follows. First, notice that the increase
in $Q^\sharp(u)$ due to the newly created waves is 
$\O(1)\cdot|\ln\ve||\Delta \Ups|$, with 
$|\Delta \Ups|=\O(1)\cdot |\sigma_{\beta_2}|$. 

Next, the difference between $Q^\sharp_{\bar{\beta}_1}$ and
$Q^\sharp_{\beta_1}$ comes from the changes in
$(W^\sharp_{\beta_1})^{-1}$ and $\tilde{z}_{\beta_1}$ which are at most 
of the order of $\sigma_{\beta_2}$ at each $x$. Hence
$$
\eqalign{
&Q^\sharp_{\bar{\beta}_1} -Q^\sharp_{\beta_1}
\le |\sigma_{\bar{\beta}_1}|\big(
\int W^\nat_{\bar{\beta}_1}(x)W^\sharp_{\bar{\beta}_1}(x) 
D_x\tilde{z}_{\bar{\beta}_1} -
\int W^\nat_{\beta_1}(x)W^\sharp_{\beta_1}(x) D_x\tilde{z}_
{\beta_1}\big) \cr
&= \O(1)\cdot |\sigma_{\beta_2}||\sigma_{\bar{\beta}_1}|
\int_0^{\tv\{u\}}{{dy}\over{(\ve +|\sigma_{\bar{\beta}_1}|+y)^2}}
+\O(1)\cdot |\sigma_{\bar{\beta}_1}|\int_0^{|\sigma_{\beta_2}|}
{{dy}\over{\ve +|\sigma_{\bar{\beta}_1}|+y}} \cr
&=\O(1)\cdot |\sigma_{\beta_2}|=
\O(1)\cdot |\Delta \Ups|.\cr
}
$$

The change in $Q^\sharp_{\alpha_i}$ also comes from the change in
$W^\sharp_{\alpha_i}(x)$ and $\tilde{z}_{\alpha_i}(x)$. 
Since the weight
$W^\sharp_{\alpha_i}(x)$ decreases as $x$  moves away from
$x_{\alpha_i}$, we have
$$\eqalign{
\sum_{\alpha_i,\, i=1,2}
|\Delta Q^\sharp_{\alpha_i}|&=\O(1)\cdot \left(
\sum_{\alpha_1}|\sigma_{\beta_2}||\sigma_{\alpha_1}|\int_{x_\beta}^\infty
(W^\sharp_{\alpha_1})^2(x) D_x\tilde z_{\alpha_1}
+\sum_{\alpha_2}|\sigma_{\beta_2}||\sigma_{\alpha_2}|\int^{x_\beta}_{-\infty}
(W^\sharp_{\alpha_2})^2(x) D_x\tilde z_{\alpha_2}
\right)\cr
&\qquad +\O(1)\cdot|\sigma_{\beta_2}|
\bigg( \sum_{\alpha_1}|\sigma_{\alpha_1}| W^\sharp_{\alpha_1}(x_{\beta_1})
+\sum_{\alpha_2} |\sigma_{\alpha_2}|W^\sharp_{\alpha_2}(x_{\beta_1})\bigg)
\cr
&=\O(1)\cdot |\ln\ve||\Delta \Ups|.\cr
}
$$

\v
\n{\bf 3.} Interaction of a shock and a rarefaction front of
different families.
To fix the ideas,
let $\beta_1$ be a shock of the $i$-th family and $\beta_2$ be
a rarefaction wave of the $j$-th family with $i>j$. Assume
$\beta_1$ and $\beta_2$ interact at time $t^*$ and denote
the outgoing  wave of the $i$-th family by $\bar{\beta}_1$,
and the $j$-th family wave $\bar{\beta}_2$. Moreover, let
$\gamma$ be a newly created shock front of the $k$-th
family, $k\ne i,j$.
By a standard interaction estimate, we have
$$
|\sigma_{\beta_i}-\sigma_{\bar{\beta}_i}|=\O(1)\cdot
|\Delta \Ups|,\quad\quad 
|\sigma_\gamma|=\O(1)\cdot|\Delta \Ups|, \quad i=1,\, 2.
$$
Thus, if we consider two shock waves $\alpha_i$, $i=1,2$
of the $k$-th family located on the left and right of the
interaction point respectively, as in the analysis of
Case 2  we have
$$\eqalign{
\sum_{\alpha_i,\, i=1,2}
|\Delta Q^\sharp_{\alpha_i}|&=\O(1)\cdot\left(
\sum_{\alpha_1}|\Delta \Ups||\sigma_{\alpha_1}|\int_{x_\beta}^\infty
(W^\sharp_{\alpha_1})^2(x) D_x\tilde z_{\alpha_1}
+\sum_{\alpha_2}|\Delta \Ups||\sigma_{\alpha_2}|\int^{x_\beta}_{-\infty}
(W^\sharp_{\alpha_2})^2(x) D_x\tilde z_{\alpha_2}
\right)\cr
&=\O(1)\cdot|\ln\ve||\Delta \Ups|.\cr
}
$$ 
In addition,
$$
| Q^\sharp_{\beta_1}-Q^\sharp_{\bar{\beta}_1}|=\O(1)\cdot
|\ln\ve||\Delta \Ups|,
\quad\quad | Q^\sharp_{\gamma}|=\O(1)\cdot|\ln\ve||\Delta \Ups|.
$$
Here we assume that $\bar{\beta}_1$ is a shock wave. In the other
case, we  have $Q^\sharp_{\beta_1}=\O(1)\cdot
|\ln\ve||\Delta \Ups|$.

\v

\n{\bf 4.} Interaction of rarefaction fronts of different families.
The change of $Q^\sharp(u)$ in this case only comes from the new
fronts 
created by the interaction. Therefore, as in the analysis of
Case 3, the total change in $Q^\sharp$ is bounded by $\O(1)\cdot|\ln\ve|
|\Delta \Ups|$.
\v

Based on the analysis of the above four cases, we see that by
choosing $C_1$ to be sufficiently large, then 
the nonlinear functional $\Hat{Q}(u)$
is non-increasing at the interaction time
when no new large shocks are introduced. 
This completes the proof of the lemma.
\endproof

\vsk
\n{\medbf 4 - Proof of the main theorem}
\v
Relying on the analysis of the two previous sections,
we can now conclude the proof of Theorem~1.
We briefly recall the main argument.  If one defines
the mollification 
$v^\delta\doteq u*\varphi_\delta$ with $\delta=\sqrt\ve$,
the estimates (1.7) hold, while (1.13)-(1.14) imply
$$\eqalign{\int_0^\tau \int &\big| v^\delta_t+A(v^\delta)
v^\delta_x-\ve v^\delta_{xx}\big|\,dxdt
=\O(1)\cdot\int_0^\tau \int\osc\big\{
u\,;~~[y-\delta,~y+\delta]\big\}\,\big|du(y)|\,dt
\cr 
&=\O(1)\cdot \int_0^\tau \sum_{|x_\alpha(t)-x_\beta(t)|\le \delta}
|\sigma_\alpha\sigma_\beta|\,dt\cr}\eqno(4.1)$$
In this case, the presence of big shocks gives a large contribution to
the right hand side (4.1), namely
$$ \int_0^\tau \sum_{\alpha\in\BS}
\big|\sigma_\alpha(t)\big|^2\,dt\eqno(4.2)$$
To get a more accurate estimate, in a neighborhood of each big shock
we replaced the mollification with a (modified) viscous travelling wave,
according to (1.21).  By doing this, we picked up more error terms, namely:
\v
\i{$\bullet$} The terms
related to the interactions of big shocks with other 
fronts.  The analysis at the beginning of Section 3 has shown
that the total contribution of all these terms satisfies the bound (1.9).
\v
\i{$\bullet$} The errors due to the difference between 
the rescaled profiles $\tilde\omega_\alpha$ in (1.19) and the exact 
travelling wave profiles $\omega_\alpha$. According to (3.4),
the total strength of these terms is
$$\int_0^\tau\sum_{\alpha\in\BS} E_\alpha(t)\,dt=\O(1)\cdot
\tau \ve\big( 1+|\ln\ve|\big)\,\tv\{\bar u\}\,.\eqno(4.3)$$
\v
On the other hand, we removed the contributions of all terms 
in (4.2).  For the function $v$ defined at (1.21) we thus have
$$\eqalign{\int_0^\tau \int \big| v_t+A(v)v_x-\ve v_{xx}\big|\,dxdt
&=
\O(1)\cdot
\tau \ve\big( 1+|\ln\ve|\big)\,\tv\{\bar u\}\cr
&\qquad +
\O(1)\cdot \int_0^\tau \bigg(\sum_{|x_\beta-x_\gamma|\le 2\sqrt\ve}
|\sigma_\beta\sigma_\gamma|-\sum_{\alpha\in\BS}|\sigma_\alpha|^2\bigg)
\,dt\,.\cr}\eqno(4.4)$$
The main goal of this section is to show that the last integral in
(4.4) can be estimated as
$$\eqalign{ \int_0^\tau \bigg(\sum_{|x_\beta-x_\gamma|\le 2\sqrt\ve}
|\sigma_\beta\sigma_\gamma|&-\sum_{\alpha\in\BS}|\sigma_\alpha|^2\bigg)
\,dt
=\O(1)\cdot\sum_{i=1}^n \int_0^\tau \bigg(
\sum_{\beta, \gamma\in\Rar_i\,,\,|x_\beta-x_\gamma|\le 8\sqrt\ve}
|\sigma_\beta\sigma_\gamma|\bigg)\,dt\cr
&\quad
+\O(1)\cdot \int_0^\tau \left|{d\over dt}\Hat Q\big(u(t)\big)\right|\,dt
+\O(1)\cdot\sqrt\ve|\ln\ve|\tau\,\tv\{\bar u\}\,.
\cr}\eqno(4.5)$$
Using the estimate (2.28) on the spreading of positive
wave-fronts and the bounds (3.15)--(3.17) concerning $\Hat Q(u)$,
from (4.5) we obtain
$$\int_0^\tau \bigg(
\sum_{|x_\beta-x_\gamma|\le 2\sqrt\ve}
|\sigma_\beta\sigma_\gamma|-\sum_{\alpha\in\BS}|\sigma_\alpha|^2\bigg)
\,dt
=\O(1)\cdot 
(1+\tau)\sqrt\ve|\ln\ve|\cdot\tv\{\bar u\}\,.$$
This will complete the proof of the estimate (1.3).
\v
The remaining part of this section is devoted to a proof of (4.5)
which is a consequence of the following lemma.

\v
\n{\bf Lemma 6.}
{\it Outside interaction times, one has ${d\over dt} \Hat Q\big(u(t))\le 0$
and  }
$$\eqalign{  \sum_{|x_\beta-x_\gamma|\le 2\sqrt\ve}
|\sigma_\beta\sigma_\gamma|-\sum_{\alpha\in\BS}|\sigma_\alpha|^2
&=\O(1)\cdot\sum_{i=1}^n  \bigg(
\sum_{\beta, \gamma\in\Rar_i\,,\,|x_\beta-x_\gamma|\le 8\sqrt\ve}
|\sigma_\beta\sigma_\gamma|\bigg)\cr
&\quad
+\O(1)\cdot  \left|{d\over dt}\Hat Q\big(u(t)\big)\right|
+\O(1)\cdot\sqrt\ve|\ln\ve|\,\tv\{\bar u\}\,.
\cr}\eqno(4.6)$$
\v
To help the reader work his way through the technicalities of the proof, 
we first describe the heart of the matter in plain words.

After removing the terms in (4.2) related to large shocks, 
the left hand side of (4.6) still contains the sum
$$\sum_{\alpha\in\SS}|\sigma_\alpha|^2,$$
where  $\SS$ denotes the set of
 all small shocks.   According to (1.22), the maximum strength
 small shock is $\leq 4\sqrt\ve\,|\ln\ve|$.
Hence the above sum is estimated by $\O(1)\cdot \sqrt\ve\,|\ln\ve|\,
\tv\{\bar u\}\,$.

Next, consider any interval $J$ of length $2\sqrt\ve$.
We first estimate the restriction of (4.6) to fronts inside $J$, i.e.
$$\Theta\doteq \sum_{x_\alpha\in J,~ |x_\alpha-x_\beta|\le 2\sqrt\ve}
|\sigma_\alpha\sigma_\beta|-\sum_{x_\alpha\in J,\,
\alpha\in\BS}|\sigma_\alpha|^2\,.
$$
It is convenient to split $\Theta$ into various sums: 
$$\eqalign{\Theta^\flat
&\doteq\sum_{x_\alpha\in J,~ |x_\alpha-x_\beta|\le 2\sqrt\ve,
~ k_\alpha\not=k_\beta}
|\sigma_\alpha\sigma_\beta|\,,\cr
\Theta_i^{\rm raref}&\doteq 
\sum_{x_\alpha\in J,~|x_\alpha-x_\beta|\le 2\sqrt\ve,
~\alpha,\beta\in \Rar_i}
|\sigma_\alpha\sigma_\beta|\,,\cr}
\qquad\qquad
\eqalign{
\Theta_i^\nat&\doteq 
\sum_{x_\alpha\in J,~ |x_\alpha-x_\beta|\le 2\sqrt\ve,
~ \alpha\in\S_i,~\beta\in \Rar_i}
|\sigma_\alpha\sigma_\beta|\,,\cr
\Theta^\sharp_i&\doteq 
\sum_{x_\alpha\in J,~|x_\alpha-x_\beta|\le 2\sqrt\ve,
~\alpha,\beta\in \S_i}
|\sigma_\alpha\sigma_\beta|\,.\cr}$$
If $\Theta^\flat$ dominates all other terms, then
the whole sum
$\Theta$ can be controlled by the rate of decrease in the functional
$Q^\flat$, related to products of fronts of different families.
The alternative case is when $J$ contains almost only waves of
one single family, say of the $i$-th family.  
If $\Theta_i^\nat$ is the dominant term, then
$\Theta$ is controlled by the decrease of the functional
$Q^\nat$ and $Q^\sharp$. If $\Theta_i^\sharp$ dominates, then $J$
contains mainly $i$-shocks, and $\Theta$ is controlled
by the decrease in $Q^\sharp$.  Finally, if $\Theta_i^{\rm raref}$ dominates,
then there is nothing to prove, because the
sum over all couples of nearby rarefactions
appears explicitly also on the right hand side of (4.6).

Covering the real line with countably many intervals
$J_\ell$ of fixed length, we eventually obtain the
desired result.
\v
\n{\bf Proof of Lemma 6.~} 
Since the system is strictly hyperbolic, the definition
of the functional $Q^\flat(u)$ implies
$$
\eqalign{
{d\over{dt}} Q^\flat(u)&=-\sum_{|x_\beta
-x_\gamma|\le 2\sqrt\ve,\, k_\beta\neq k_\gamma} |\sigma_\beta\sigma_\gamma|
{{|\dot{x}_\beta-\dot{x}_\gamma|}\over {4\sqrt\ve}}
\cr
&=-{c_2\over\sqrt\ve} \sum_{|x_\beta-x_\gamma|\le 2\sqrt\ve,
\, k_\beta\neq k_\gamma} |\sigma_\beta\sigma_\gamma|,
\cr
}\eqno(4.7)
$$
for some constant $c_2>0$ related to the minimum gap
between different characteristic speeds.  Hence the terms 
containing a product of two waves of different families 
on the left hand side of (4.6) are
controlled by the decreasing rate of $Q^\flat(u)$. 

In the following, we only need to show that 
the products involving one or two
shock waves of the same family can be controlled by
the decreasing rate of the nonlinear functional $\Hat{Q}(u)$,
plus the quantity in
(2.28) and $\sqrt\ve|\ln\ve|\tv\{\bar{u}\}$.

By the definition of $Q^\nat(u)$, we know that the rarefaction
waves located in $I_\alpha(t)$ involved
in $Q^\nat_\alpha(u)$ approach to the large shock
wave $\alpha$ unless there are waves of the other families in
between. Hence, if we use $\BS'$ to denote the set of big
shocks $\alpha$ such that the total strength
of small wave fronts within the interval $I_\alpha(t)$ is 
$\le {|\sigma_\alpha|
\over 4}$, then for $\alpha\in \BS'$, we have
$$
\eqalign{
 {d\over {dt}}Q_\alpha^\nat(u)&=-\sum_{x_\beta\in I_\alpha(t),\,
\beta\in \Rar_\alpha} |\sigma_\beta|{|\dot{x}_\alpha-
\dot{x}_\beta|\over {4\sqrt\ve}} \cr
&\le -{c_3\over\sqrt\ve}\sum_{x_\beta\in
I_\alpha(t),\,\beta\in \Rar_\alpha} |\sigma_\alpha\sigma_\beta|
+{\O(1)\over\sqrt\ve}\sum_{\beta, x_\beta\in
I_\alpha(t),\, k_\beta\neq k_\alpha} |\sigma_\alpha\sigma_\beta|,
\cr
}
\eqno(4.8)
$$
where
$$
Q^\nat_\alpha(u)=\int W^\nat_\alpha(x) D_x\tilde{w}_\alpha.
$$

On the other hand, the functional 
$Q^\sharp(u)$ is defined for all shock waves no matter they are
small or large. In this way, its time derivative yields
mainly the product of two shock waves of the same family
with distance $\leq 2\sqrt\ve$. Let
$$
Q_\alpha^\sharp(u)
\doteq |\sigma_\alpha|\,
\int W_\alpha^\nat(x)\,W_\alpha^\sharp(x) D_x\tilde z_\alpha\,.
$$
Since there is a factor $3$ in front of the summation of
rarefaction waves in the definition of $Q^\sharp(u)$,
this guarantees that all the shock waves appearing in
$Q^\sharp_\alpha$ approach to the  shock wave $\alpha$
if there is no waves of other families in between.
Notice that there is a constant $\ve$ in the denominator
of the weight function $W^\sharp_\alpha(x)$.
Thus, for any shock $\alpha\in \S$,
$$
\eqalign{
{d\over {dt}} Q_\alpha^\sharp(u)
&\le
-{c_4\over \sqrt\ve}
\Hat{\sum}_{\beta \in \S_\alpha,\,x_\beta\in \hat{I}_\alpha(t)}
|\sigma_\alpha\sigma_\beta|+
{c_5\over\sqrt\ve}
\sum_{\beta \in \Rar_\alpha,\,x_\beta\in I_\alpha(t)}
|\sigma_\alpha\sigma_\beta|
+{\O(1)\over\sqrt\ve}\sum_{\beta, x_\beta\in
I_\alpha(t),\, k_\beta\neq k_\alpha} |\sigma_\alpha\sigma_\beta|,
\cr
&\le
-{c_4\over \sqrt\ve}
\sum_{\beta \in \S_\alpha,\,x_\beta\in \hat{I}_\alpha(t)}
|\sigma_\alpha\sigma_\beta|+
{c_5\over\sqrt\ve}
\sum_{\beta \in \Rar_\alpha,\,x_\beta\in I_\alpha(t)}
|\sigma_\alpha\sigma_\beta|\cr
&\qquad
+{\O(1)\over\sqrt\ve}\sum_{\beta, x_\beta\in
I_\alpha(t),\, k_\beta\neq k_\alpha} |\sigma_\alpha\sigma_\beta|
+\O(1)\sqrt\ve |\sigma_\alpha|,
\cr
}
\eqno(4.9)
$$
where $c_4, c_5>0$ are constants independent of $\ve$,
$\hat{I}_\alpha(t)=[x_\alpha-2\sqrt\ve, x_\alpha[\,\bigcup\,
]x_\alpha, x_\alpha+2\sqrt\ve]$,
and $\Hat{\sum}$ means that the  summation is over all shocks $\beta$ 
with the property that the total strength of all shock fronts between
$\alpha$ and $\beta$ with $x_\beta\in I_\alpha(t)$
is  $\geq\ve$.

By noticing that the time derivative of $\Ups(u)$ is zero outside
interaction times and by choosing $C_2>\!>C_3>\!> 1$,
based on the estimates (4.7)-(4.9), the increase of
(4.6) can be given as follows, by considering separately
the products involving
large shocks, and those involving only small wave fronts.

For a large shock front $\alpha$, consider the summation
$$
\sum_{x_\beta\in \hat{I}_\alpha(t)}|\sigma_\alpha\sigma_\beta|.
\eqno(4.10)
$$
Since 
the sum of all products  $|\sigma_\alpha\sigma_\beta|$
when $k_\alpha\neq k_\beta$ is controlled by (4.7), we have 
$$
\sum_{
x_\beta\in \hat{I}_\alpha(t)}|\sigma_\alpha\sigma_\beta|
\le \O(1)\cdot\left|{d\over dt}Q^\flat_\alpha(u)\right|\,,
\eqno(4.11)
$$
provided that
waves of different families dominate, say 
$$
\sum_{\beta, x_\beta\in I_\alpha(t), k_\beta\neq k_\alpha}
|\sigma_\beta|\ge {1\over 4} 
\sum_{\beta, x_\beta\in \hat{I}_\alpha(t)}|\sigma_\beta|.
\eqno(4.12)
$$

It thus remains to consider the case when (4.12) does not hold.
We then have
$$
\sum_{x_\beta\in \hat{I}_\alpha(t), k_\beta =k_\alpha} |\sigma_\beta|
\ge {3\over 4} 
\sum_{\beta, x_\beta\in \hat{I}_\alpha(t)}|\sigma_\beta|.
\eqno(4.13)
$$
In this case, 
if $\alpha\in \BS'$, then the summation of
$|\sigma_\alpha\sigma_\beta|$ for $\beta\in \Rar_\alpha\cup\S_\alpha$
is controlled by (4.8) and (4.9) together with (4.7).
Therefore
$$
\eqalign{
\sum_{\alpha\in \BS'}\,
\sum_{x_\beta\in \hat{I}_\alpha(t)}|\sigma_\alpha\sigma_\beta|
&\le \O(1)\cdot\bigg( \sum_{i=1}^n 
\sum_{\beta, \gamma\in\Rar_i\,,\,|x_\beta-x_\gamma|\le 2\sqrt\ve}
|\sigma_\beta\sigma_\gamma|\bigg)\cr
&\qquad
+\O(1)\cdot \left|{d\over dt}\Hat Q\big(u(t)\big)\right|
+\O(1)\cdot \ve \,\tv\{\bar u\}\,.
\cr}\eqno(4.14)
$$
Moreover, by (4.9), for $\alpha\in \S$, if 
$$
\sum_{\beta\in \S_\alpha,\, x_\beta\in \hat{I}_\alpha(t)}
|\sigma_\beta|\ge {{2c_5}\over c_4}
\sum_{\beta\in \Rar_\alpha,\, x_\beta\in \hat{I}_\alpha(t)}
|\sigma_\beta|,
\eqno(4.15)
$$
then
$$
\sum_{
x_\beta\in \hat{I}_\alpha(t)}|\sigma_\alpha\sigma_\beta|
\le \O(1)\cdot
\bigg(-{d\over dt}\big( C_2\sqrt\ve|\ln\ve|Q^\nat_\alpha 
+\sqrt\ve Q^\sharp_\alpha\big)
+\ve |\sigma_\alpha| \bigg)\,.
\eqno(4.16)
$$

For a large shock wave, it now remains to consider the case when
$\alpha\in \BS''=\BS-\BS'$ satisfying (4.13) and
$$
\sum_{\beta\in \S_\alpha, x_\beta\in \hat{I}_\alpha(t)}
|\sigma_\beta|\le {{2c_5}\over c_4}
\sum_{\beta\in \Rar_\alpha,\, x_\beta\in \hat{I}_\alpha(t)}
|\sigma_\beta|.
$$
We denote the set consisting all these large shock waves by $\BS'''$.
Notice that this is a subset of $\BS''$.
Roughly speaking, for  $\alpha\in \BS'''$, the small wave fronts is
not small compared to $\alpha$  and rarefaction waves of $k_\alpha$-th
family dominate in $I_\alpha(t)$.
Hence, for $\alpha\in \BS'''$, one has
$$
\sum_{\theta\in \BS'''_\alpha, x_\theta\in I_\alpha(t)}\,\,
\sum_{\beta, x_\beta\in \hat{I}_\theta(t)}
|\sigma_\theta\sigma_\beta|
\le \O(1)\cdot\sum_{\beta,\gamma\in \Rar_\alpha, \,x_\beta, x_\gamma\in 
[x_\alpha -4\sqrt\ve, x_\alpha +4\sqrt\ve]}
|\sigma_\beta\sigma_\gamma|\,,
\eqno(4.17)
$$
which is controlled by the corresponding part of (2.28) in the
interval $[x_\alpha-4\sqrt\ve, ~x_\alpha +4\sqrt\ve]$. Hence
$$
\sum_{\alpha\in\BS'''}\,\, \sum_{\beta,
x_\beta\in \hat{I}_\alpha(t)}
|\sigma_\alpha\sigma_\beta|
\le \O(1)\cdot\sum_{\alpha,\beta\in\Rar, \,|x_\alpha-x_\beta|\le 8\sqrt\ve,
 k_\alpha=k_\beta}
|\sigma_\alpha\sigma_\beta|\,,
\eqno(4.18)
$$
which is estimated by (2.28).

Combining (4.11), (4.14),  (4.16) and (4.18) we obtain

$$
\eqalign{
\sum_{\alpha\in \BS}\,
\sum_{x_\beta\in \hat{I}_\alpha(t)}|\sigma_\alpha\sigma_\beta|
&\le \O(1)\cdot\bigg( \sum_{i=1}^n 
\sum_{\beta, \gamma\in\Rar_i\,,\,|x_\beta-x_\gamma|\le 8\sqrt\ve}
|\sigma_\beta\sigma_\gamma|\bigg)\cr
&\qquad
+\O(1)\cdot \left|{d\over dt}\Hat Q\big(u(t)\big)\right|
+\O(1)\cdot \ve \,\tv\{\bar u\}\,.
\cr}\eqno(4.19)
$$

Now it remains to show the sum
of products of small wave fronts of the same
family satisfies the same bound:
$$
\eqalign{
\sum_{\alpha,\beta\in\SS\cup \Rar,
|x_\alpha -x_\beta|\le 2\sqrt\ve,\, k_\alpha=k_\beta}
|\sigma_\alpha\sigma_\beta|
&\le \O(1)\cdot\bigg( \sum_{i=1}^n 
\sum_{\beta, \gamma\in\Rar_i\,,\,|x_\beta-x_\gamma|\le 8\sqrt\ve}
|\sigma_\beta\sigma_\gamma|\cr
&\quad
+ \left|{d\over dt}\Hat Q\big(u(t)\big)\right|
+\sqrt\ve|\ln\ve|\,\tv\{\bar u\}\,\bigg).\cr}
\eqno(4.20)
$$
To obtain the estimate (4.20), we divide the real line into
a union of closed 
 intervals of length $2\sqrt\ve$, i.e.
$\R=\bigcup_{i} J_i$ with $J_i\doteq \big[2i\sqrt\ve,\,2(i+1)
\sqrt\ve\big]$.
We denote by $s_i^k$ and $r_i^k$ 
respectively the total strengths of small $k$-shock  and $k$-rarefaction
fronts contained in the interval
in $J_i$.
We have
$$
\sum_{\alpha,\beta\in\SS\cup\Rar,
|x_\alpha -x_\beta|\le 2\sqrt\ve,\, k_\alpha=k_\beta}
|\sigma_\alpha\sigma_\beta|
\le \sum_i\sum_{\alpha\in \SS\cup\Rar,
x_\alpha\in J_i}\,\,\sum_{\beta\in\SS\cup\Rar,
|x_\alpha-x_\beta|\le 2\sqrt\ve,
k_\alpha=k_\beta}
|\sigma_\alpha\sigma_\beta|.
$$
To estimate the quantity
$$
\sum_{\alpha\in \SS\cup\Rar, 
x_\alpha\in J_i}\,\,\sum_{\beta\in\SS\cup\Rar,
|x_\alpha-x_\beta|\le 2\sqrt\ve,
k_\alpha=k_\beta}
|\sigma_\alpha\sigma_\beta|,
\eqno(4.21)
$$
we consider the following two cases.

\v

\n$\bullet\quad$ For a given $k$,
 $s^k_i\ge {{2 c_5}\over c_4}(r^k_{i-1}+r^k_i +r^k_{i+1})$.
In this case,
from (4.9) we deduce 
$$
\eqalign{
&{d\over {dt}} \sum_{\alpha\in \SS_k,\, x_\alpha\in J_i} Q^\sharp_\alpha(u)
\le
-{c_5\over \sqrt\ve}\sum_{\alpha\in\SS_k,\, x_\alpha\in J_i}\,\,
\sum_{\beta \in \S_\alpha,\,x_\beta\in \hat{I}_\alpha(t)}
|\sigma_\alpha\sigma_\beta|\cr
&\qquad +{\O(1)\over\sqrt\ve}\sum_{\alpha\in\SS_k,\, x_\alpha\in J_i}\,\,
\sum_{x_\beta\in
I_\alpha(t),\, k_\beta\neq k_\alpha} |\sigma_\alpha\sigma_\beta|
+\O(1)\sqrt\ve \sum_{\alpha\in \SS_k, x_\alpha\in J_i} |\sigma_\alpha|\cr
&\le
-{c_5\over \sqrt\ve}\sum_{\alpha\in\SS_k,\, x_\alpha\in J_i}\,\,
\sum_{\beta \in \S_\alpha,\,x_\beta\in I_\alpha(t)}
|\sigma_\alpha\sigma_\beta|
\qquad +{\O(1)\over\sqrt\ve}\sum_{\alpha\in\SS_k,\, x_\alpha\in J_i}\,\,
\sum_{x_\beta\in
I_\alpha(t),\, k_\beta\neq k_\alpha} |\sigma_\alpha\sigma_\beta|\cr
&\qquad
+{\O(1)\over \sqrt\ve}\sum_{\alpha\in SS_k, x_\alpha\in J_i}|\sigma_\alpha|^2
+\O(1)\sqrt\ve \sum_{\alpha\in \SS_k, x_\alpha\in J_i} |\sigma_\alpha|\cr
&\le
-{c_5\over \sqrt\ve}\sum_{\alpha\in\SS_k,\, x_\alpha\in J_i}\,\,
\sum_{\beta \in \S_\alpha,\,x_\beta\in I_\alpha(t)}
|\sigma_\alpha\sigma_\beta|
+{\O(1)\over\sqrt\ve}\sum_{\alpha\in\SS_k,\, x_\alpha\in J_i}\,\,
\sum_{x_\beta\in
I_\alpha(t),\, k_\beta\neq k_\alpha} |\sigma_\alpha\sigma_\beta|\cr
&\qquad +\O(1)|\ln\ve|\sum_{\alpha\in SS_k, x_\alpha\in J_i}|\sigma_\alpha|
 +\O(1)\sqrt\ve \sum_{\alpha\in \SS_k, x_\alpha\in J_i} |\sigma_\alpha|.\cr}
\eqno(4.22)
$$
Here we have used the fact that $\sigma_\alpha\le 4\sqrt\ve|\ln\ve|$
for $\alpha\in \SS_k$. 
By (4.22) we see that those terms containing a  product
with $\alpha\in\SS_k$
and $\beta\in\SS_k$ in (4.21) can be controlled by ${d\over dt} \Hat{Q}(u)$
up to an error of the order of $\sqrt\ve|\ln\ve|\tv\{\bar{u}\}$.
Since the total stength of all small $k$-shocks in $J_i$ dominates
the total strength of all $k$-rarefactions 
in $\bigcup_{j=i-1}^{i+1} J_j$, the products
of  $\alpha\in \SS_k$ and $\beta\in\Rar_k$, and the products
of $\alpha\in\Rar_k$ with $\beta\in\SS_k$ for $ x_\beta\in J_i$ in (4.21),
are also  controlled by  ${d\over dt} \Hat{Q}(u)$
up to an error of the order of $\sqrt\ve|\ln\ve|\tv\{\bar{u}\}$.
 Moreover, those
products of $\alpha\in \Rar_k$ and $\beta\in \Rar_k$ in (4.21) are
 controlled by the corresponding parts of (2.28) the interval
$\bigcup_{j=i-1}^{i+1} J_i$.

Hence, it remains to consider the
product of $\alpha\in\Rar_k$ and 
$\beta\in\SS_k$ with $ x_\beta\in J_{i-1}\cup
J_{i+1}$. 
To fix the ideas, we consider the case when $\alpha\in \Rar_k$,
$\beta\in \SS_k$ with $x_\beta\in J_{i-1}$, i.e.,
$$
\sum_{\alpha\in \Rar_k,\beta\in \SS_k,\, x_\alpha\in J_{i},\,
x_\beta\in J_{i-1},\, |x_\alpha-x_\beta|\le 2\sqrt\ve} |\sigma_\alpha\sigma_\beta|.
\eqno(4.23)
$$
When $s^k_{i-1}
\le {{2 c_5}\over c_4}(r^k_{i-2}+r^k_{i-1} +r^k_{i})$, (4.23) is controlled
by (2.28) in the interval $\bigcup_{j=i-2}^{i+1} J_j$. 
Otherwise
$$
\sum_{\alpha\in \Rar_k,\beta\in \SS_k,\, x_\alpha\in J_{i},\,
x_\beta\in J_{i-1},\, |x_\alpha-x_\beta|\le 2\sqrt\ve} 
|\sigma_\alpha\sigma_\beta|\le \bigg(\sum_{\beta,
x_\beta\in J_{i-1}, \beta\in\SS_k}|\sigma_\beta|\bigg)^2,$$
which can be controlled as in (4.22), using 
(4.7) to control the $k$-shock
fronts in $J_{i-1}$.

\v

\n$\bullet\quad$ Now assume that 
$s^k_i< {{2 c_5}\over c_4}(r^k_{i-1}+r^k_i +r^k_{i+1})$.
In this case the total strength of all $k$-rarefaction fronts 
in $\bigcup_{j=i-1}^{i+1} J_j$
dominates the total strength of $k$-small shocks in $J_i$.
As done previously, we only need to
consider the case when $\alpha\in\SS_k\cup\Rar_k$ and $\beta\in\SS_k$ with $
x_\beta\in J_{i-1}\cup J_{i+1}$ in (4.21) because 
all the other terms can be
controlled by (2.28) in the corresponding interval 
$\bigcup_{j=i-1}^{i+1} J_j$.

For illustration, we discuss the following two terms,
$$\sum_{\alpha\in \Rar_k,\beta\in \SS_k,\, x_\alpha\in J_{i},\,
x_\beta\in J_{i-1},\, |x_\alpha-x_\beta|\le 2\sqrt\ve} 
|\sigma_\alpha\sigma_\beta|,
\eqno(4.24)
$$
and
$$\sum_{\alpha,\beta\in \SS_k,\, x_\alpha\in J_{i},\,
x_\beta\in J_{i-1},\, |x_\alpha-x_\beta|\le 2\sqrt\ve} 
|\sigma_\alpha\sigma_\beta|,
\eqno(4.25)
$$
respectively as follows. The other terms can be handled 
similarly.

Concerning (4.24), when  $s^k_{i-1}
\le {{2 c_5}\over c_4}(r^k_{i-2}+r^k_{i-1} +r^k_{i})$,
it can be controlled by the corresponding terms in (2.28)
in the interval $\bigcup_{j=i-2}^{i} J_j$. Otherwise,
when $s^k_{i-1}
> {{2 c_5}\over c_4}(r^k_{i-2}+r^k_{i-1} +r^k_{i})$,
we have
$$
\sum_{\alpha\in \Rar_k,\,\beta\in \SS_k,\, x_\alpha\in J_{i},\,
x_\beta\in J_{i-1},\, |x_\alpha-x_\beta|\le 2\sqrt\ve} 
|\sigma_\alpha\sigma_\beta|
\le
\bigg(\sum_{\beta\in \SS_k,\, x_\beta\in J_{i-1}} |\sigma_\beta|\bigg)^2,
$$
which can be estimated as in (4.22), using 
(4.7) to control the $k$-shock fronts in the interval $J_{i-1}$.

Concerning (4.25), when  $s^k_{i-1}
\ge {{2 c_5}\over c_4}(r^k_{i-2}+r^k_{i-1} +r^k_{i})$,
then a similar argument as in (4.22) can be applied, using
(4.7).

Otherwise,
$$
\sum_{\alpha,\beta\in \SS_k,\, x_\alpha\in J_{i},\,
x_\beta\in J_{i-1},\, |x_\alpha-x_\beta|
\le 2\sqrt\ve} |\sigma_\alpha\sigma_\beta|
\le \bigg(\sum_{j=i-2}^{i+1} r^k_j\bigg)^2,
\eqno(4.26)
$$
which can be controlled by the corresponding term in (2.28)
in the interval $\bigcup_{j=i-2}^{i+1} J_j$.

Notice that each interval $J_i$ can be counted no more than three times.
By combining (4.22)-(4.26), we have desired estimate on (4.21) for small
wave fronts so that (4.20) holds.
In summary, (4.19) and (4.20) imply 
(4.6), completing
the  proof of the lemma.

\endproof
\v
\n{\bf Remark 4.} In the proof of the error estimate
(1.3), the three basic ingredients are:

\i{$\bullet$} The existence of uniformly Lipschitz semigroups
of approximate (viscous) solutions.

\i{$\bullet$} The decay of positive waves, due to genuine nonlinearity,

\i{$\bullet$} The exponential rate of convergence to
steady states, in the tails of travelling viscous shocks.

\n Assuming that all characteristic fields are genuinely nonlinear,
we thus conjecture that similar error estimates are valid also
for the semidiscrete scheme considered in [Bi].
In the case of straight line systems, based on the analysis
in [BJ], 
it is reasonable to expect that analogous results
should also hold
for the Godunov scheme.
\v
\n{\bf Remark 5.}  In the case where all characteristic fields
are linearly degenerate, 
solutions with Lipschitz continuous initial data
having small total variation remain uniformly
Lipschitz continuous for all times, as shown in [B1].
Therefore, the easy error estimate 
(1.16) can be used.   For systems having some linenearly
degenerate and some genuinely nonlinear fields,
we still conjecture that the error bound (1.3)
is valid.  A proof, however, will require some new techniques.
Indeed, the contact discontinuities that may be generated
by shock interactions at times $t>0$ can no longer be
approximated by viscous travelling profiles.

\vsk

{\bf Acknowledgments:} 
The first author was supported by the Italian M.I.U.R.,
within the research project \# 2002017219  ``Equazioni iperboliche e
paraboliche non lineari''. The research of the second author was  supported
by the Competitive Earmarked Research Grant of Hong Kong
CityU 1142/01P \# 9040648.

\vsk

\c{\medbf References} \v
\i{[Bi]} S.~Bianchini,
BV solutions for the semidiscrete upwind scheme, {\it Arch. Rat. Mech. 
Anal.} {\bf 167} (2003), 1-81.
\v
\i{[BiB]} S.~Bianchini and A.~Bressan, Vanishing viscosity
solutions of nonlinear hyperbolic systems, {\it Annals of
Mathematics}, to appear. 
\v
\i{[B1]} A.~Bressan,
Contractive metrics for nonlinear hyperbolic systems, {\it Indiana Univ.
J. Math.} {\bf 37} (1988), 409-421.
\v
\i{[B2]} A.~Bressan,
{\it Hyperbolic Systems of Conservation Laws. The One Dimensional
Cauchy Problem}, Oxford University Press, 2000.
\v
\i{[BC1]} A.~Bressan and R.~M.~Colombo, 
The semigroup generated by $2\times 2$ conservation laws,
{\it Arch. Rat. Mech. Anal.} {\bf 113} (1995), 1-75 
\v
\i{[BC2]} A.~Bressan and R.~M.~Colombo, Decay of positive waves in
nonlinear systems of conservation laws, {\it Ann. Scuola Norm.
Sup. Pisa} {\bf IV - 26} (1998), 133-160. 
\v
\i{[BCP]} A.~Bressan, G.~Crasta and B.~Piccoli, Well posedness of
the Cauchy problem for $n\times n$ conservation laws, {\it Amer.
Math. Soc. Memoir} {\bf 694} (2000).
\v
\i{[BLY]} A.~Bressan, T.~P.~Liu and T.~Yang,  $ L^1$ stability
estimates for $n\times n$ conservation laws, {\it Arch. Rational
Mech. Anal.} {\bf 149} (1999), 1-22. 
\v
\i{[BM]} A.~Bressan and A.~Marson, Error bounds for a
deterministic version of the Glimm scheme, {\it Arch. Rat. Mech.
Anal.} {\bf 142} (1998), 155-176. 
\v
\i{[BJ]} A.~Bressan and H.~K.~Jenssen,
On the convergence of Godunov scheme for nonlinear hyperbolic
systems,
{\it Chinese Ann. Math.} {\bf B - 21} (2000), 1-16.
\v
\i{[BY]} A.~Bressan and T.~Yang, A sharp decay estimate for
positive nonlinear waves, {\it SIAM J. Math. Anal.}, submitted.
\v
\i{[GX]} J.~Goodman and Z.~Xin, Viscous limits for piecewise
smooth solutions to systems of conservation laws, {\it Arch.
Rational Mech. Anal.} {\bf 121} (1992), 235-265.
\v
\i{[K]} N.~N.~Kuznetsov, Accuracy of some approximate methods for computing
the weak solutions of a first-order quasi-linear equation, 
{\it U.S.S.R. Comp. Math. and Math. Phys.} {\bf 16} (1976), 105-119.
\v
\i{[L]} T.~P.~Liu, Admissible solutions of hyperbolic conservation
laws, {\it Amer. Math. Soc. Memoir} {\bf 240} (1981). 
\v
\i{[O]} O.~Oleinik, Discontinuous solutions of nonlinear
differential equations, {\it Amer. Math. Soc. Transl.} {\bf 26}
(1963), 95-172. 
\v
\i{[TT]} T.~Tiang and Z.~H.~Teng, The sharpness of Kuznetsov's 
$\O(\sqrt{\Delta x})$  ~$L^1$ error estimate for monotone difference
schemes, {\it Math. Comp.} {\bf 64} (1995), 581-589.

\bye